\documentstyle[12pt]{article}        
\title{A classification of inner actions of the 
Dipper-Donkin quantization $GL_2$ on the 
Clifford algebra ${\it C}(1,3)$.}  
\author{      
Suemi Rodr\'\i guez-Romo\\Centre of Theoretical Research\\
National University of M\'exico, Campus Cuautitl\'an\\
Apdo. Postal 95, Unidad Militar, Cuautitl\'an Izcalli\\
Estado de M\'exico, 54768 M\'exico.$^{\ast}$}         
\date{}  
\begin{document}   
\maketitle
\renewcommand{\thefootnote}{\fnsymbol{footnote}}
\setcounter{footnote}{-1}
\footnote{$\hspace*{-6mm}^{\ast}$
e-mail: suemi@servidor.unam.mx}
\renewcommand{\thefootnote}{\arabic{footnote}}
\baselineskip0.6cm
{\small {\bf Abstract.} We present all inner actions on the 
Clifford algebra ${\it C}(1,3)$ of the quantum group $GL_2$ 
constructed by Dipper and Donkin \cite{dido}.\vspace{2cm}

\baselineskip0.6cm  
Following the method already developed for studying the actions of 
$GL_q(2,C)$ on the Clifford algebra ${\it C}(1,3)$ \cite{aqg} we 
construct any possible inner action of Dipper-Donkin quantization 
on this Clifford algebra. We also provide with the corresponding 
operator algebra $\Re$ (namely the image of the representation),
the algebra of invariants $I$ which is equal to the centralizer of
$\Re$ in ${\it C}(1,3)$, and the {\it perturbation} of the representation
\cite{suemi}. Let $\varphi$ be a representation and $c_{11}$, $c_{12}$, 
$c_{21}$, $c_{22}$, generators of $GL_2$, we say that the representation 
$\varphi$ has a {\it perturbation} if $\varphi(c_{12}c_{21})\neq 0$. 
We define $c_{ij}\rightarrow C_{ij}$ to be a finite dimensional 
representation of the quantum $GL_2$.\

\centerline{\Huge  CASE 1) $d$=$diag(q^2,q,1,1)$}
\vspace*{1.5cm}
\vbox{\offinterlineskip
\hrule
\halign{&\vrule#&\strut \quad \hfil# \hfil
 &\vrule#&\quad \hfil# \hfil&\vrule#& \quad# \hfil \cr
height5pt&\omit&& \omit&& \omit&\cr&
        {\bf CASE 1.1)}   
&&
$\begin{array}{l}
C_{12}=\alpha e_{12}+\beta e_{23}+\gamma e_{2
4}\\ 
C_{21}=0 \\
C_{11}={\bf 1} \\
C_{22}=q^2e_{11}+qe_{22}+e_{33}+e_{44}\\
\end{array}$
&& 
        $\Re =\left(\matrix{*&*&*&*\cr 0&*&*&*\cr 0&0&
\epsilon &0\cr 0&0&0&\epsilon \cr}
\right)$&\cr
height5pt&\omit&& \omit&&\omit&\cr
\noalign{\hrule }
height5pt&\omit&& \omit&&\omit&\cr    
&$\matrix{dim \Re \cr 8 \cr 
dim I \cr 1}$ &&
         $\begin{array}{l}
I\cong C \end{array}$ && 
         Perturbation zero &\cr
height5pt&\omit&& \omit&& \omit&\cr
\noalign{\hrule }
height2pt&\omit&& \omit&& \omit&\cr
\noalign{\hrule }
height5pt&\omit&& \omit&& \omit&\cr&
        {\bf CASE 1.2)}               
&&
$\begin{array}{c}
C_{12}=\alpha e_{12}+\beta e_{24}\\
C_{21}=0\\
C_{11}={\bf 1}+e_{34}\\
C_{22}=q^2e_{11}+qe_{22}+e_{33}+e_{44}-e_{34}
\end{array}
$
&& 
        $\Re =\left(\matrix{*&*&0&0\cr 0&*&0&*\cr 0&0&
\epsilon &*\cr 0&0&0&\epsilon \cr}
\right)$&\cr
height5pt&\omit&& \omit&&\omit&\cr
\noalign{\hrule }
height5pt&\omit&& \omit&&\omit&\cr    
&
$\matrix{dim \Re \cr 6 \cr dim I \cr 2}$
 &&
$I =\left(\matrix{\alpha&0&0&0\cr 0&\alpha&0&0\cr 0&0&
\alpha & \gamma\cr 0&0&0&\alpha \cr}
\right)$
&& 
         Perturbation zero &\cr
height5pt&\omit&& \omit&& \omit&\cr
\noalign{\hrule }
height2pt&\omit&& \omit&& \omit&\cr
\noalign{\hrule }
height5pt&\omit&& \omit&& \omit&\cr&
{\bf CASE 1.3)}   
&&
$\begin{array}{c}
	C_{12}=0\\
	C_{21}=\alpha e_{21}+\beta e_{32}\\
	C_{11}=e_{11}+q^{-1}e_{22}+q^{-2}e_{33}+\gamma e_{44}\\
        C_{22}=q^{2}e_{11}+q^{2}e_{22}+q^{2}e_{33}+\gamma^{-1}e_{44}
\end{array}$
&& 
        $\Re =\left(\matrix{*&0&0&0\cr *&*&0&0\cr 0&*&
\epsilon &0\cr 0&0&0&\epsilon \cr}
\right)$&\cr
height5pt&\omit&& \omit&&\omit&\cr
\noalign{\hrule }
height5pt&\omit&& \omit&&\omit&\cr    
&$\matrix{dim \Re \cr 5 \cr 
dim I \cr 3}$ 
         &&
$I =\left(\matrix{\alpha&0&0&0\cr 0&\alpha&0&0\cr 0&0&
\alpha & \beta\cr 0&0&0&\alpha \cr}
\right)$
&& 
         Perturbation zero&\cr
height5pt&\omit&& \omit&& \omit&\cr
\noalign{\hrule }
height2pt&\omit&& \omit&& \omit&\cr
\noalign{\hrule }}}
\vbox{\offinterlineskip
\hrule
\halign{&\vrule#&\strut \quad \hfil# \hfil
 &\vrule#&\quad \hfil# \hfil&\vrule#& \quad# \hfil \cr
height5pt&\omit&& \omit&& \omit&\cr
        &{\bf CASE 1.4)}  
&&
$\begin{array}{l}
C_{12}=0\\
C_{21}=\alpha e_{21}+\beta e_{32}\\
C_{11}=e_{11}+q^{-1}e_{22}+q^{-2}e_{33}+q^{-2}e_{44}+e_{34}\\ 
C_{22}=q^{2}{\bf 1}-q^4e_{34}\end{array}
$
&& 
        $\Re =\left(\matrix{*&0&0&0\cr *&*&0&0\cr 0&*&
\epsilon &*\cr 0&0&0&\epsilon \cr}
\right)$&\cr
height5pt&\omit&& \omit&&\omit&\cr
\noalign{\hrule }
height5pt&\omit&& \omit&&\omit&\cr    
&$\matrix{dim \Re \cr 6 \cr 
dim I \cr 2}$ && 
$I =\left(\matrix{\alpha&0&0&0\cr 0&\alpha&0&0\cr 0&0&
\alpha & \beta\cr 0&0&0&\alpha \cr}
\right)$&& 
Perturbation zero &\cr
height5pt&\omit&& \omit&& \omit&\cr
\noalign{\hrule }
height2pt&\omit&& \omit&& \omit&\cr
\noalign{\hrule }
height5pt&\omit&& \omit&& \omit&\cr
        & {\bf CASE 1.5)}   
&&
$\begin{array}{c}
	C_{12}=\beta e_{23}+\gamma e_{24}\\
	C_{21}=\alpha e_{21}\\
	C_{11}=diag(1, q^{-1}, q^{-1}, q^{-1})\\
        C_{22}=diag(q^{2}, q^{2}, q, q)
\end{array}$
&& 
        $\Re =\left(\matrix{*&0&0&0\cr *
&*&\varphi&-\frac{\beta\varphi}{\gamma}\cr 0&0&
\epsilon &0\cr 0&0&0&\epsilon \cr}
\right)$&\cr
height5pt&\omit&& \omit&&\omit&\cr
\noalign{\hrule }
height5pt&\omit&& \omit&&\omit&\cr    
&$\matrix{dim \Re \cr 5 \cr 
dim I \cr 1}$ && 
$I\cong C$
           && 
Perturbation zero &\cr
height5pt&\omit&& \omit&& \omit&\cr
\noalign{\hrule }
height2pt&\omit&& \omit&& \omit&\cr
\noalign{\hrule }
height5pt&\omit&& \omit&& \omit&\cr
        & {\bf CASE 1.6)}   
&&
$\begin{array}{l}
C_{12}=\alpha e_{24} \\ 
C_{21}=\beta e_{21} \\ 
C_{11}= e_{11}+q^{-1}e_{22}+q^{-1}e_{33}+q^{-1}e_{44}+e_{34}\\
C_{22}= q^2e_{11}+q^{2}e_{22}+qe_{33}+qe_{44}-q^2e_{34}
\end{array}
$
&& 
        $\Re =\left(\matrix{*&0&0&0\cr *&*&0&*\cr 0&0&
\epsilon &*\cr 0&0&0&\epsilon \cr}
\right)$&\cr
height5pt&\omit&& \omit&&\omit&\cr
\noalign{\hrule }
height5pt&\omit&& \omit&&\omit&\cr    
&$\matrix{dim \Re \cr 6 \cr 
dim I \cr 3}$ && 
$I =\left(\matrix{\alpha&0&0&0\cr 0&\alpha&0&0\cr 0&0&
\alpha & \beta\cr 0&0&\gamma&\alpha \cr}
\right)$
           && 
Perturbation zero &\cr
height5pt&\omit&& \omit&& \omit&\cr
\noalign{\hrule }
height2pt&\omit&& \omit&& \omit&\cr
\noalign{\hrule }
height5pt&\omit&& \omit&& \omit&\cr
        & {\bf CASE 1.7)}   
&&
$\begin{array}{l}
C_{12}=\alpha e_{12} \\ 
C_{21}=\beta e_{32}+\gamma e_{42} \\ 
C_{11}= e_{11}+e_{22}+q^{-1}e_{33}+q^{-1}e_{44}\\
C_{22}= q^2e_{11}+qe_{22}+qe_{33}+qe_{44}
\end{array}
$
&& 
        $\Re =\left(\matrix{*&*&0&0\cr 0&*&0&0\cr 0&*&
\epsilon &0\cr 0&*&0&\epsilon \cr}
\right)$&\cr
height5pt&\omit&& \omit&&\omit&\cr
\noalign{\hrule }
height5pt&\omit&& \omit&&\omit&\cr    
&$\matrix{dim \Re \cr 6 \cr 
dim I \cr 1}$ && 
$I\cong C$
           && 
Perturbation zero &\cr
height5pt&\omit&& \omit&& \omit&\cr
\noalign{\hrule }
height2pt&\omit&& \omit&& \omit&\cr
\noalign{\hrule }}}
\vbox{\offinterlineskip
\hrule
\halign{&\vrule#&\strut \quad \hfil# \hfil
 &\vrule#&\quad \hfil# \hfil&\vrule#& \quad# \hfil \cr
height5pt&\omit&& \omit&& \omit&\cr
	&{\bf CASE 1.8)}   
&&
$\begin{array}{l}
C_{12}=\alpha e_{12} \\ 
C_{21}=\beta e_{32}\\ 
C_{11}= e_{11}+e_{22}+q^{-1}e_{33}+q^{-1}e_{44}+e_{34}\\
C_{22}= q^2e_{11}+qe_{22}+qe_{33}+qe_{44}-q^2e_{34}
\end{array}
$
&& 
        $\Re =\left(\matrix{*&*&0&0\cr 0&*&0&0\cr 0&*&
\epsilon &*\cr 0&0&0&\epsilon \cr}
\right)$&\cr
height5pt&\omit&& \omit&&\omit&\cr
\noalign{\hrule }
height5pt&\omit&& \omit&&\omit&\cr    
&$\matrix{dim \Re \cr 6 \cr 
dim I \cr 2 \cr }$ && 
$I =\left(\matrix{\alpha&0&0&0\cr 0&\alpha&0&0\cr 0&0&
\alpha & \beta\cr 0&0&0&\alpha \cr}
\right)$
           && 
Perturbation zero &\cr
height5pt&\omit&& \omit&& \omit&\cr
\noalign{\hrule }}}


\centerline{\Huge  CASE 2) $d$=$diag(q^2,q,q,1)$}
\vspace*{1.5cm}
\vbox{\offinterlineskip
\hrule
\halign{&\vrule#&\strut \quad \hfil# \hfil
 &\vrule#&\quad \hfil# \hfil&\vrule#& \quad# \hfil \cr
height5pt&\omit&& \omit&& \omit&\cr&
        {\bf CASE 2.1)}   
&&
$\begin{array}{l}
C_{12}=0 \\
C_{21}=-\mu\delta e_{21}+\mu\gamma e_{31}+\gamma e_{42}+\delta e_{43}\\ 
C_{11}=e_{11}+q^{-1}e_{22}+q^{-1}e_{33}+q^{-2}e_{44} \\
C_{22}=q^2{\bf 1}\\
\end{array}$
&& 
        $\Re =\left(\matrix{*&0&0&0\cr \varphi&\epsilon&0&0\cr 
\frac{\gamma\varphi}{\delta}&0&0&0\cr 
0&\phi&\frac{-\gamma\phi}{\delta}&* \cr}
\right)$&\cr
height5pt&\omit&& \omit&&\omit&\cr
\noalign{\hrule }
height5pt&\omit&& \omit&&\omit&\cr    
&$\matrix{dim \Re \cr 5 \cr 
dim I \cr 5}$ &&
         $I =\left(\matrix{\beta&0&0&0\cr 0&\psi&\phi&0\cr 
0&\frac{\gamma}{\delta}(\varphi-\psi)+
\left(\frac{\gamma}{\delta}\right)^2\phi&\varphi &0\cr 0&0&0&\alpha \cr}
\right)$
&& 
         Perturbation zero &\cr
height5pt&\omit&& \omit&& \omit&\cr
\noalign{\hrule }
height2pt&\omit&& \omit&& \omit&\cr
\noalign{\hrule }
height5pt&\omit&& \omit&& \omit&\cr&
        {\bf CASE 2.2)}               
&&
$\begin{array}{c}
C_{12}=qe_{13}-\mu e_{24}\\
C_{21}=-\mu e_{21}+e_{43}\\
C_{11}=e_{11}+q^{-1}e_{22}+e_{33}+q^{-1}e_{44}\\
C_{22}=q^2e_{11}+q^2e_{22}+qe_{33}+qe_{44}-q\mu e_{23}
\end{array}
$
&& 
        $\Re =\left(\matrix{*&*\cr 0&*\cr}\right)\otimes 
\left(\matrix{*&*\cr 0&*\cr}\right)$&\cr
height5pt&\omit&& \omit&&\omit&\cr
\noalign{\hrule }
height5pt&\omit&& \omit&&\omit&\cr    
&
$\matrix{dim \Re \cr 9 \cr dim I \cr 1}$
 &&
$I\cong C$
&& 

$\matrix{{\rm Perturbation} \cr 
C_{12}C_{21}=-\mu e_{23}}$&\cr
height5pt&\omit&& \omit&& \omit&\cr
\noalign{\hrule }
height2pt&\omit&& \omit&& \omit&\cr
\noalign{\hrule }
height5pt&\omit&& \omit&& \omit&\cr&
{\bf CASE 2.3)}   
&&
$\begin{array}{c}
	C_{12}=q\lambda\delta e_{13}\\
	C_{21}=\delta e_{43}\\
	C_{11}=e_{11}+\alpha e_{22}+e_{33}+q^{-1}e_{44} \\
        C_{22}=q^{2}e_{11}+q\alpha^{-1}e_{22}+qe_{33}+qe_{44} \\
	\alpha\neq 1 \\
\end{array}$
&& 
        $\Re =\left(\matrix{*&0&*&0\cr 0&*&0&0\cr 0&0&*&0\cr 
0&0&*&*\cr}
\right)$&\cr
height5pt&\omit&& \omit&&\omit&\cr
\noalign{\hrule }
height5pt&\omit&& \omit&&\omit&\cr    
&$\matrix{dim \Re \cr 6 \cr 
dim I \cr 2}$ 
         &&
$I =\left(\matrix{\alpha&0&0&0\cr 0&\beta&0&0\cr 0&0&
\alpha & 0\cr 0&0&0&\alpha \cr}
\right)$
&& 
         Perturbation zero&\cr
height5pt&\omit&& \omit&& \omit&\cr
\noalign{\hrule }
height2pt&\omit&& \omit&& \omit&\cr
\noalign{\hrule }}}
\vbox{\offinterlineskip
\hrule
\halign{&\vrule#&\strut \quad \hfil# \hfil
&\vrule#&\quad \hfil# \hfil&\vrule#& \quad# \hfil \cr
height5pt&\omit&& \omit&& \omit&\cr
&{\bf CASE 2.4)}  
&&
$\begin{array}{l}
C_{12}=q\lambda\delta e_{13}\\
C_{21}=\delta e_{43}\\
C_{11}=e_{11}+e_{22}+e_{33}+q^{-1}e_{44}  \\ 
C_{22}=q^2e_{11}+qe_{22}+qe_{33}+qe_{44}\\
\end{array}$
&& 
        $\Re =\left(\matrix{*&0&*&0\cr 0&\epsilon&0&0\cr 0&0&
\epsilon &0\cr 0&0&*&* \cr}
\right)$&\cr
height5pt&\omit&& \omit&&\omit&\cr
\noalign{\hrule }
height5pt&\omit&& \omit&&\omit&\cr    
&$\matrix{dim \Re \cr 5 \cr 
dim I \cr 2}$ && 
$I =\left(\matrix{\alpha&0&0&0\cr 0&\beta&0&0\cr 0&0&
\alpha &0\cr 0&0&0&\alpha \cr}
\right)$&& 
Perturbation zero &\cr
height5pt&\omit&& \omit&& \omit&\cr
\noalign{\hrule }
height2pt&\omit&& \omit&& \omit&\cr
\noalign{\hrule }
height5pt&\omit&& \omit&& \omit&\cr
        & {\bf CASE 2.5)}   
&&
$\begin{array}{c}
	C_{12}=q\gamma e_{12}+q\delta e_{13}\\
	C_{21}=\gamma e_{42}+\delta e_{43}\\
	C_{11}=e_{11}+e_{22}+e_{33}+q^{-1}e_{44}\\
        C_{22}=q^2e_{11}+qe_{22}+q_{33}+qe_{44}
\end{array}$
&& 
        $\Re =\left(\matrix{*&\varphi&\varphi&0\cr 0
&\epsilon&0&0\cr 0&0&
\epsilon &0\cr 0&\psi&\psi&*\cr}
\right)$&\cr
height5pt&\omit&& \omit&&\omit&\cr
\noalign{\hrule }
height5pt&\omit&& \omit&&\omit&\cr    
&$\matrix{dim \Re \cr 5 \cr 
dim I \cr 6}$ && 
$I =\left(\matrix{\alpha&0&0&0\cr 0&\gamma&\delta&0\cr 0&\phi&
\psi &0\cr 0&0&0&\beta\cr}
\right)$
           && 
Perturbation zero &\cr
height5pt&\omit&& \omit&& \omit&\cr
\noalign{\hrule }
height2pt&\omit&& \omit&& \omit&\cr
\noalign{\hrule }
height5pt&\omit&& \omit&& \omit&\cr
        & {\bf CASE 2.6)}   
&&
$\begin{array}{l}
C_{12}=0\\ 
C_{21}=0\\ 
C_{11}= e_{11}+\alpha e_{22}+\beta e_{33}+\gamma e_{44}\\
C_{22}= q^2e_{11}+q\alpha^{-1}e_{22}+
q\beta^{-1}e_{33}+\gamma^{-1}e_{44}\\
\alpha\neq \beta\neq \gamma\neq 1\end{array}
$
&& 
        $\Re =\left(\matrix{*&0&0&0\cr 0&*&0&0\cr 0&0&
*&0\cr 0&0&0&* \cr}
\right)$&\cr
height5pt&\omit&& \omit&&\omit&\cr
\noalign{\hrule }
height5pt&\omit&& \omit&&\omit&\cr    
&$\matrix{dim \Re \cr 4 \cr 
dim I \cr 4}$ && 
$I =\left(\matrix{\alpha&0&0&0\cr 0&\beta&0&0\cr 0&0&
\gamma &0\cr 0&0&0&\delta \cr}
\right)$
           && 
Perturbation zero &\cr
height5pt&\omit&& \omit&& \omit&\cr
\noalign{\hrule }
height2pt&\omit&& \omit&& \omit&\cr
\noalign{\hrule }
height5pt&\omit&& \omit&& \omit&\cr

        & {\bf CASE 2.7)}   
&&
$\begin{array}{l}
C_{12}=0\\ 
C_{21}=0 \\ 
C_{11}= e_{11}+\alpha e_{22}+\alpha e_{33}+\gamma e_{44}\\
C_{22}= q^2e_{11}+q\alpha^{-1}e_{22}+q\alpha^{-1}e_{33}
+\gamma^{-1}e_{44}
\end{array}
$
&& 
$\Re =\left(\matrix{*&0&0&0\cr 0&\epsilon&0&0\cr 0&0&
\epsilon &0\cr 0&0&0&*\cr}
\right)$
        &\cr
height5pt&\omit&& \omit&&\omit&\cr
\noalign{\hrule }
height5pt&\omit&& \omit&&\omit&\cr    
&$\matrix{dim \Re \cr 3 \cr 
dim I \cr 6}$ && 
$I =\left(\matrix{\alpha&0&0&0\cr 0&\gamma&\delta&0\cr 0&\epsilon&
\phi & 0\cr 0&0&0&\beta \cr}
\right)$
           && 
Perturbation zero &\cr
height5pt&\omit&& \omit&& \omit&\cr
\noalign{\hrule }
height2pt&\omit&& \omit&& \omit&\cr
\noalign{\hrule }}}
\vbox{\offinterlineskip
\hrule
\halign{&\vrule#&\strut \quad \hfil# \hfil
 &\vrule#&\quad \hfil# \hfil&\vrule#& \quad# \hfil \cr
height5pt&\omit&& \omit&& \omit&\cr

	&{\bf CASE 2.8)}   
&&
$\begin{array}{l}
C_{12}=qe_{12}+\mu e_{34} \\ 
C_{21}=\mu e_{31}+e_{42}\\ 
C_{11}= e_{11}+e_{22}+q^{-1}e_{33}+q^{-1}e_{44}\\
C_{22}= q^2e_{11}+qe_{22}+q^2e_{33}+qe_{44}+q\mu e_{32}
\end{array}
$
&&
$\Re =\left(\matrix{*&*\cr 0&*\cr}
\right)\otimes \left(\matrix{*&0\cr *&*\cr}
\right)$
        &\cr
height5pt&\omit&& \omit&&\omit&\cr
\noalign{\hrule }
height5pt&\omit&& \omit&&\omit&\cr    
&$\matrix{dim \Re \cr 9 \cr 
dim I \cr 1 \cr }$ && 
$I \cong C$
           && 
$\matrix{{\rm Perturbation} \cr 
C_{12}C_{21}=\mu e_{32}}$&\cr
height5pt&\omit&& \omit&& \omit&\cr
\noalign{\hrule }
height2pt&\omit&& \omit&& \omit&\cr
\noalign{\hrule }
height5pt&\omit&& \omit&& \omit&\cr

	&{\bf CASE 2.9)}   
&&
$\begin{array}{l}
C_{12}=0\\ 
C_{21}=0\\ 
C_{11}= e_{11}+\alpha e_{22}+\alpha e_{33}+\beta e_{44}\\
C_{22}= q^2e_{11}+\frac{q}{\alpha}e_{22}+\frac{q}{\alpha}e_{33}+
\beta^{-1}e_{44}
\end{array}
$
&&
$\Re =\left(\matrix{*&0&0& 0\cr 0&\epsilon&0&0\cr
0&0&\epsilon&0\cr0&0&0&*\cr}\right)$
        &\cr
height5pt&\omit&& \omit&&\omit&\cr
\noalign{\hrule }
height5pt&\omit&& \omit&&\omit&\cr    
&$\matrix{dim \Re \cr 3\cr
dim I \cr 6 \cr }$ && 
$I =\left(\matrix{\alpha&0&0&0\cr 0&\gamma&\delta&0\cr 0&\phi&
\psi & 0\cr 0&0&0&\beta \cr}
\right)$
           && 
Perturbation zero &\cr
height5pt&\omit&& \omit&& \omit&\cr
\noalign{\hrule }
height2pt&\omit&& \omit&& \omit&\cr
\noalign{\hrule }
height5pt&\omit&& \omit&& \omit&\cr

	&{\bf CASE 2.10)}   
&&
$\begin{array}{l}
C_{12}=0\\ 
C_{21}=0\\ 
C_{11}= e_{11}+\alpha e_{22}+\beta e_{33}+\gamma e_{44}\\
C_{22}= q^2e_{11}+\frac{q}{\alpha}e_{22}+\frac{q}{\beta}e_{33}+
\gamma^{-1}e_{44}\\
		\\
\alpha\neq \beta
\end{array}
$
&&
$\Re =\left(\matrix{*&0&0& 0\cr 0&*&0&0\cr
0&0&*&0\cr0&0&0&*\cr}\right)$
        &\cr
height5pt&\omit&& \omit&&\omit&\cr
\noalign{\hrule }
height5pt&\omit&& \omit&&\omit&\cr    
&$\matrix{dim \Re \cr 4 \cr 
dim I \cr 4 \cr }$ && 
$I =\left(\matrix{\alpha&0&0&0\cr 0&\beta&0&0\cr 0&0&
\gamma & 0\cr 0&0&0&\delta \cr}\right)$
           && 
Perturbation zero &\cr
height5pt&\omit&& \omit&& \omit&\cr
height5pt&\omit&& \omit&& \omit&\cr
\noalign{\hrule }
height2pt&\omit&& \omit&& \omit&\cr
\noalign{\hrule }
height5pt&\omit&& \omit&& \omit&\cr

	&{\bf CASE 2.11)}   
&&
$\begin{array}{l}
C_{12}=qe_{12}+\mu e_{34}\\ 
C_{21}=\mu e_{31}+e_{42}\\ 
C_{11}= e_{11}+e_{22}+q^{-1}e_{33}+q^{-1}e_{44}\\
C_{22}= q^2e_{11}+qe_{22}+q^2e_{33}+qe_{44}+q\mu e_{32}
\end{array}
$
&&
$\Re =\left(\matrix{*&*\cr 0&*\cr}\right)\otimes
\left(\matrix{*&0\cr *&*\cr}\right)$
        &\cr
height5pt&\omit&& \omit&&\omit&\cr
\noalign{\hrule }
height5pt&\omit&& \omit&&\omit&\cr    
&$\matrix{dim \Re \cr 9\cr
dim I \cr 1 \cr }$ && 
$I \cong C$
           && 
$\matrix{{\rm Perturbation}\cr
C_{12}C_{21}=\mu e_{32}\cr }$&\cr
height5pt&\omit&& \omit&& \omit&\cr
\noalign{\hrule }
height2pt&\omit&& \omit&& \omit&\cr
\noalign{\hrule }}}
\vbox{\offinterlineskip
\hrule
\halign{&\vrule#&\strut \quad \hfil# \hfil
 &\vrule#&\quad \hfil# \hfil&\vrule#& \quad# \hfil \cr
height5pt&\omit&& \omit&& \omit&\cr

	&{\bf CASE 2.12)}   
&&
$\begin{array}{l}
C_{12}=0 \\ 
C_{21}=-\mu\delta e_{21}+\delta e_{43}\\ 
C_{11}= e_{11}+q^{-1}e_{22}+q^{-1}e_{33}+q^{-2}e_{44}+e_{23}\\
C_{22}= q^2e_{11}+q^2e_{22}+q^2e_{33}+q^2e_{44}-q^3e_{23}
\end{array}
$
&&
$\Re =\left(\matrix{*&0&0&0\cr *&\epsilon&*&0\cr
0&0&\epsilon&0\cr0&0&*&*\cr}\right)$
        &\cr
height5pt&\omit&& \omit&&\omit&\cr
\noalign{\hrule }
height2pt&\omit&& \omit&&\omit&\cr    
&$\matrix{dim \Re \cr 6\cr
dim I \cr 1 \cr }$ && 
$I \cong C$
           && 
Perturbation zero &\cr
height5pt&\omit&& \omit&& \omit&\cr
\noalign{\hrule }
height2pt&\omit&& \omit&& \omit&\cr
\noalign{\hrule }
height5pt&\omit&& \omit&& \omit&\cr

	&{\bf CASE 2.13)}   
&&
$\begin{array}{l}
C_{12}=0\\ 
C_{21}=0\\ 
C_{11}= e_{11}+\alpha e_{22}+\alpha e_{33}+\beta e_{44}+e_{23}\\
C_{22}= q^2e_{11}+\frac{q}{\alpha}e_{22}+\frac{q}{\alpha}e_{33}+
\frac{1}{\beta}e_{44}-\frac{q}{\alpha^2}e_{23}
\end{array}
$
&&
$\Re =\left(\matrix{*&0&0& 0\cr 0&\epsilon&*&0\cr
0&0&\epsilon&0\cr0&0&0&*\cr}\right)$
        &\cr
height5pt&\omit&& \omit&&\omit&\cr
\noalign{\hrule }
height5pt&\omit&& \omit&&\omit&\cr    
&$\matrix{dim \Re \cr 4\cr
dim I \cr 4 \cr }$ && 
$I =\left(\matrix{\alpha&0&0&0\cr 0&\gamma&\delta&0\cr 0&0&
\gamma & 0\cr 0&0&0&\beta \cr}\right)$
           && 
Perturbation zero &\cr
height5pt&\omit&& \omit&& \omit&\cr
\noalign{\hrule }
height2pt&\omit&& \omit&& \omit&\cr
\noalign{\hrule }
height5pt&\omit&& \omit&& \omit&\cr

	&{\bf CASE 2.14)}   
&&
$\begin{array}{l}
C_{12}=q\lambda\delta e_{13}\\ 
C_{21}=\delta e_{43}\\ 
C_{11}= e_{11}+e_{22}+e_{33}+q^{-1}e_{44}+e_{23}\\
C_{22}= q^2e_{11}+qe_{22}+qe_{33}+qe_{44}-qe_{23}
\end{array}
$
&&
$\Re =\left(\matrix{*&0&*&0\cr 0&\epsilon&*&0\cr
0&0&\epsilon&0\cr 0&0&*&*\cr}\right)$
        &\cr
height5pt&\omit&& \omit&&\omit&\cr
\noalign{\hrule }
height5pt&\omit&& \omit&&\omit&\cr    
&$\matrix{dim \Re \cr 6\cr
dim I \cr 2 \cr }$ && 
$I =\left(\matrix{\alpha&0&0&0\cr 0&\alpha&\beta&0\cr 0&0&
\alpha& 0\cr 0&0&0&\alpha \cr}\right)$
           && 
Perturbation zero &\cr
height5pt&\omit&& \omit&& \omit&\cr
\noalign{\hrule }
height2pt&\omit&& \omit&& \omit&\cr
\noalign{\hrule }
height5pt&\omit&& \omit&& \omit&\cr

	&{\bf CASE 2.15)}   
&&
$\begin{array}{l}
C_{12}=0\\ 
C_{21}=0\\ 
C_{11}= e_{11}+\alpha e_{22}+\alpha e_{33}+\beta e_{44}+e_{23}\\
C_{22}= q^2e_{11}+\frac{q}{\alpha}e_{22}+\frac{q}{\alpha}e_{33}+
\frac{1}{\beta}e_{44}-\frac{q}{\alpha^2}e_{23}
\end{array}
$
&&
$\Re =\left(\matrix{*&0&0&0\cr 0&\epsilon&*&0\cr
0&0&\epsilon&0\cr0&0&0&*\cr}\right)$
        &\cr
height5pt&\omit&& \omit&&\omit&\cr
\noalign{\hrule }
height5pt&\omit&& \omit&&\omit&\cr    
&$\matrix{dim \Re \cr 4\cr
dim I \cr 3 \cr }$ && 
$I =\left(\matrix{\alpha&0&0&0\cr 0&\beta&\gamma&0\cr 0&0&
\beta& 0\cr 0&0&0&\alpha \cr}\right)$
           && 
Perturbation zero &\cr
height5pt&\omit&& \omit&& \omit&\cr
\noalign{\hrule }}}


\centerline{\Huge  CASE 3) $d$=$diag(q^2,q^2,q,1)$}
\vspace*{1.5cm}
\vbox{\offinterlineskip
\hrule
\halign{&\vrule#&\strut \quad \hfil# \hfil
 &\vrule#&\quad \hfil# \hfil&\vrule#& \quad# \hfil \cr
height5pt&\omit&& \omit&& \omit&\cr&
        {\bf CASE 3.1)}   
&&
$\begin{array}{l}
C_{12}=\alpha e_{13}+\beta e_{23}+\gamma e_{34} \\
C_{21}=0\\ 
C_{11}={\bf 1}\\
C_{22}=q^2e_{11}+q^2e_{22}+qe_{33}+e_{44}\\
\end{array}$
&& 
        $\Re =\left(\matrix{\epsilon&0&\phi&\varphi\cr 
0&\epsilon&\phi&\varphi\cr 
0&0&*&*\cr 
0&0&0&* \cr}
\right)$&\cr
height5pt&\omit&& \omit&&\omit&\cr
\noalign{\hrule }
height5pt&\omit&& \omit&&\omit&\cr    
&$\matrix{dim \Re \cr 6 \cr 
dim I \cr 5}$ &&
         $I =\left(\matrix{\varphi&\beta&0&0\cr 
\gamma&\delta&0&0\cr 0&0&\alpha&0\cr 0&0&0&\alpha \cr}
\right)$
&& 
         Perturbation zero &\cr
height5pt&\omit&& \omit&& \omit&\cr
\noalign{\hrule }
height2pt&\omit&& \omit&& \omit&\cr
\noalign{\hrule }
height5pt&\omit&& \omit&& \omit&\cr&
        {\bf CASE 3.2)}               
&&
$\begin{array}{c}
C_{12}=\alpha e_{13}+\gamma e_{34}\\
C_{21}=0\\
C_{11}={\bf 1}+e_{12}\\
C_{22}=q^2e_{11}+q^2e_{22}+e_{33}+e_{44}
\end{array}
$
&& 
        $\Re =\left(\matrix{\epsilon&*&*&*\cr 
0&\epsilon&0&0\cr 0&0&*&*\cr 0&0&0&*\cr}\right)$&\cr
height5pt&\omit&& \omit&&\omit&\cr
\noalign{\hrule }
height5pt&\omit&& \omit&&\omit&\cr    
&
$\matrix{dim \Re \cr 7 \cr dim I \cr 2}$
 &&
$I =\left(\matrix{\alpha&\beta&0&0\cr 
0&\alpha&0&0\cr 0&0&\alpha&0\cr 0&0&0&\alpha \cr}
\right)$
&& Perturbation zero&\cr
height5pt&\omit&& \omit&& \omit&\cr
\noalign{\hrule }
height2pt&\omit&& \omit&& \omit&\cr
\noalign{\hrule }
height5pt&\omit&& \omit&& \omit&\cr&
{\bf CASE 3.3)}   
&&
$\begin{array}{c}
	C_{12}=0\\
	C_{21}=\alpha e_{31}+\beta e_{32}+\gamma e_{43}\\
	C_{11}=q^{2}e_{11}+q^{2}e_{22}+qe_{33}+e_{44}\\
        C_{22}={\bf 1}
\end{array}$
&& 
        $\Re =\left(\matrix{\epsilon&0&0&0\cr 0&\epsilon&0&0\cr 
\varphi&\varphi&*&0\cr 
\eta&\eta&*&*\cr}
\right)$&\cr
height5pt&\omit&& \omit&&\omit&\cr
\noalign{\hrule }
height5pt&\omit&& \omit&&\omit&\cr    
&$\matrix{dim \Re \cr 6 \cr 
dim I \cr 5}$ 
         &&
$I =\left(\matrix{\beta&\gamma&0&0\cr \delta&\eta&0&0\cr 0&0&
\alpha & 0\cr 0&0&0&\alpha \cr}
\right)$
&& 
         Perturbation zero&\cr
height5pt&\omit&& \omit&& \omit&\cr
\noalign{\hrule }
height2pt&\omit&& \omit&& \omit&\cr
\noalign{\hrule }}}
\vbox{\offinterlineskip
\hrule
\halign{&\vrule#&\strut \quad \hfil# \hfil
&\vrule#&\quad \hfil# \hfil&\vrule#& \quad# \hfil \cr
height5pt&\omit&& \omit&& \omit&\cr
&{\bf CASE 3.4)}  
&&
$\begin{array}{l}
C_{12}=0\\
C_{21}=\beta e_{32}+\gamma e_{43}\\
C_{11}={\bf 1}+e_{12}\\ 
C_{22}=q^2e_{11}+q^{2}e_{22}+qe_{33}+e_{44}-q^2e_{12}\end{array}
$
&& 
        $\Re =\left(\matrix{\epsilon&*&0&0\cr 0&\epsilon&0&0\cr 0&*&
* &0\cr 0&*&*&* \cr}
\right)$&\cr
height5pt&\omit&& \omit&&\omit&\cr
\noalign{\hrule }
height5pt&\omit&& \omit&&\omit&\cr    
&$\matrix{dim \Re \cr 7 \cr 
dim I \cr 2}$ && 
$I =\left(\matrix{\beta&0&0&0\cr 0&\alpha&0&0\cr 0&0&
\alpha &0\cr 0&0&0&\alpha \cr}
\right)$&& 
Perturbation zero &\cr
height5pt&\omit&& \omit&& \omit&\cr
\noalign{\hrule }
height2pt&\omit&& \omit&& \omit&\cr
\noalign{\hrule }
height5pt&\omit&& \omit&& \omit&\cr
        & {\bf CASE 3.5)}   
&&
$\begin{array}{c}
	C_{12}=\alpha e_{13}+\beta e_{23}\\
	C_{21}=\gamma e_{43}\\
	C_{11}=qe_{11}+qe_{22}+qe_{33}+e_{44}\\
        C_{22}=qe_{11}+qe_{22}+e_{33}+e_{44}
\end{array}$
&& 
        $\Re =\left(\matrix{\epsilon&0&\varphi&0\cr 0
&\epsilon&\varphi&0\cr 0&0&
*&0\cr 0&0&*&*\cr}
\right)$&\cr
height5pt&\omit&& \omit&&\omit&\cr
\noalign{\hrule }
height5pt&\omit&& \omit&&\omit&\cr    
&$\matrix{dim \Re \cr 5 \cr 
dim I \cr 5}$ && 
$I =\left(\matrix{\beta&\gamma&0&0\cr \epsilon&\varphi&0&0\cr 
0&0&\alpha&0\cr 0&0&0&\alpha\cr}
\right)$
           && 
Perturbation zero &\cr
height5pt&\omit&& \omit&& \omit&\cr
\noalign{\hrule }
height2pt&\omit&& \omit&& \omit&\cr
\noalign{\hrule }
height5pt&\omit&& \omit&& \omit&\cr
        & {\bf CASE 3.6)}   
&&
$\begin{array}{l}
C_{12}=\alpha e_{13}+\beta e_{23}\\ 
C_{21}=\gamma e_{43}\\ 
C_{11}= qe_{11}+qe_{22}+qe_{33}+e_{44}+e_{12}\\
C_{22}= qe_{11}+qe_{22}+e_{33}+e_{44}-qe_{12}
\end{array}
$
&& 
        $\Re =\left(\matrix{\epsilon&*&*&0\cr 0&\epsilon&*&0\cr 0&0&
*&0\cr 0&0&0&* \cr}
\right)$&\cr
height5pt&\omit&& \omit&&\omit&\cr
\noalign{\hrule }
height5pt&\omit&& \omit&&\omit&\cr    
&$\matrix{dim \Re \cr 6 \cr 
dim I \cr 1}$ && 
$I\cong C$
           && 
Perturbation zero &\cr
height5pt&\omit&& \omit&& \omit&\cr
\noalign{\hrule }
height2pt&\omit&& \omit&& \omit&\cr
\noalign{\hrule }
height5pt&\omit&& \omit&& \omit&\cr
        & {\bf CASE 3.7)}   
&&
$\begin{array}{l}
C_{12}=\gamma e_{34}\\ 
C_{21}=\beta e_{32} \\ 
C_{11}= qe_{11}+qe_{22}+e_{33}+e_{44}\\
C_{22}= qe_{11}+qe_{22}+qe_{33}+e_{44}
\end{array}
$
&& 
$\Re =\left(\matrix{\epsilon&0&0&0\cr 0&\epsilon&0&0\cr 0&*&
*&*\cr 0&0&0&*\cr}
\right)$
        &\cr
height5pt&\omit&& \omit&&\omit&\cr
\noalign{\hrule }
height5pt&\omit&& \omit&&\omit&\cr    
&$\matrix{dim \Re \cr 5 \cr 
dim I \cr 3}$ && 
$I =\left(\matrix{\beta&\gamma&0&0\cr 0&\alpha&0&0\cr 0&0&
\alpha& 0\cr 0&0&0&\alpha \cr}
\right)$
           && 
Perturbation zero &\cr
height5pt&\omit&& \omit&& \omit&\cr
\noalign{\hrule }
height2pt&\omit&& \omit&& \omit&\cr
\noalign{\hrule }}}
\vbox{\offinterlineskip
\hrule
\halign{&\vrule#&\strut \quad \hfil# \hfil
 &\vrule#&\quad \hfil# \hfil&\vrule#& \quad# \hfil \cr
height5pt&\omit&& \omit&& \omit&\cr
	&{\bf CASE 3.8)}   
&&
$\begin{array}{l}
C_{12}=\gamma e_{34} \\ 
C_{21}=\beta e_{32}\\ 
C_{11}= \alpha e_{11}+qe_{22}+e_{33}+e_{44}\\
C_{22}= \frac{q^2}{\alpha}e_{11}+qe_{22}+qe_{33}+e_{44}\\
\alpha\neq q
\end{array}
$
&&
$\Re =\left(\matrix{*&0&0&0\cr 
0&*&0&0\cr 0&*&*&*\cr 0&0&0&*\cr}
\right)$
        &\cr
height5pt&\omit&& \omit&&\omit&\cr
\noalign{\hrule }
height5pt&\omit&& \omit&&\omit&\cr    
&$\matrix{dim \Re \cr 6\cr 
dim I \cr 1 \cr }$ && 
$I \cong C$
           && 
Perturbation zero&\cr
height5pt&\omit&& \omit&& \omit&\cr
\noalign{\hrule }
height2pt&\omit&& \omit&& \omit&\cr
\noalign{\hrule }
height5pt&\omit&& \omit&& \omit&\cr

	&{\bf CASE 3.9)}   
&&
$\begin{array}{l}
C_{12}=\gamma e_{34}\\ 
C_{21}=\beta e_{32}\\ 
C_{11}= qe_{11}+qe_{22}+e_{33}+e_{44}+e_{12}\\
C_{22}= qe_{11}+qe_{22}+qe_{33}+e_{44}-e_{12}
\end{array}
$
&&
$\Re =\left(\matrix{\epsilon&*&0& 0\cr 0&\epsilon&0&0\cr
0&*&*&*\cr0&0&0&*\cr}\right)$
        &\cr
height5pt&\omit&& \omit&&\omit&\cr
\noalign{\hrule }
height5pt&\omit&& \omit&&\omit&\cr    
&$\matrix{dim \Re \cr 6\cr
dim I \cr 2 \cr }$ && 
$I =\left(\matrix{\alpha&\beta&0&0\cr 0&\alpha&0&0\cr 0&0&
\alpha& 0\cr 0&0&0&\alpha \cr}
\right)$
           && 
Perturbation zero &\cr
height5pt&\omit&& \omit&& \omit&\cr
\noalign{\hrule }
height2pt&\omit&& \omit&& \omit&\cr
\noalign{\hrule }
height5pt&\omit&& \omit&& \omit&\cr

	&{\bf CASE 3.10)}   
&&
$\begin{array}{l}
C_{12}=\gamma e_{34}\\ 
C_{21}=\alpha e_{31}\\ 
C_{11}= qe_{11}+qe_{22}+e_{33}+e_{44}\\
C_{22}= qe_{11}+qe_{22}+qe_{33}+e_{44}
\end{array}
$
&&
$\Re =\left(\matrix{\epsilon&0&0& 0\cr 0&\epsilon&0&0\cr
*&0&*&*\cr0&0&0&*\cr}\right)$
        &\cr
height5pt&\omit&& \omit&&\omit&\cr
\noalign{\hrule }
height5pt&\omit&& \omit&&\omit&\cr    
&$\matrix{dim \Re \cr 5\cr
dim I \cr 3 \cr }$ && 
$I =\left(\matrix{\alpha&0&0&0\cr \beta&\gamma&0&0\cr 0&0&
\alpha & 0\cr 0&0&0&\alpha \cr}\right)$
           && 
Perturbation zero &\cr
height5pt&\omit&& \omit&& \omit&\cr
height5pt&\omit&& \omit&& \omit&\cr
\noalign{\hrule }
height2pt&\omit&& \omit&& \omit&\cr
\noalign{\hrule }
height5pt&\omit&& \omit&& \omit&\cr

	&{\bf CASE 3.11)}   
&&
$\begin{array}{l}
C_{12}=\gamma e_{34}\\ 
C_{21}=\beta e_{31}\\ 
C_{11}= qe_{11}+\alpha e_{22}+e_{33}+e_{44}\\
C_{22}= qe_{11}+\frac{q^2}{\alpha}e_{22}+qe_{33}+e_{44}\\
\alpha\neq q
\end{array}
$
&&
$\Re =\left(\matrix{*&0&0&0\cr 0&*&0&0\cr
*&0&*&*\cr 0&0&0&*\cr}\right)$
        &\cr
height5pt&\omit&& \omit&&\omit&\cr
\noalign{\hrule }
height5pt&\omit&& \omit&&\omit&\cr    
&$\matrix{dim \Re \cr 6\cr
dim I \cr 1 \cr }$ && 
$I \cong C$
           && 
Perturbation zero &\cr
height5pt&\omit&& \omit&& \omit&\cr
\noalign{\hrule }
height2pt&\omit&& \omit&& \omit&\cr
\noalign{\hrule }}}
\vbox{\offinterlineskip
\hrule
\halign{&\vrule#&\strut \quad \hfil# \hfil
 &\vrule#&\quad \hfil# \hfil&\vrule#& \quad# \hfil \cr
height5pt&\omit&& \omit&& \omit&\cr

	&{\bf CASE 3.12)}   
&&
$\begin{array}{l}
C_{12}=\gamma e_{34} \\ 
C_{21}=0 \\ 
C_{11}= \alpha e_{11}+\alpha e_{22}+e_{33}+e_{44}+e_{12}\\
C_{22}= \frac{q^2}{\alpha}e_{11}+\frac{q^2}{\alpha}e_{22}+qe_{33}+
e_{44}-\frac{q^2}{\alpha^2}e_{12}
\end{array}
$
&&
$\Re =\left(\matrix{\epsilon&*&0&0\cr 0&\epsilon&0&0\cr
0&0&*&*\cr0&0&0&*\cr}\right)$
        &\cr
height5pt&\omit&& \omit&&\omit&\cr
\noalign{\hrule }
height2pt&\omit&& \omit&&\omit&\cr    
&$\matrix{dim \Re \cr 5\cr
dim I \cr 3 \cr }$ && 
$I =\left(\matrix{\beta&\gamma&0&0\cr 0&\beta&0&0\cr 0&0&
\alpha & 0\cr 0&0&0&\alpha \cr}\right)$
           && 
Perturbation zero &\cr
height5pt&\omit&& \omit&& \omit&\cr
\noalign{\hrule }
height2pt&\omit&& \omit&& \omit&\cr
\noalign{\hrule }
height5pt&\omit&& \omit&& \omit&\cr

	&{\bf CASE 3.13)}   
&&
$\begin{array}{l}
C_{12}=\gamma e_{34}\\ 
C_{21}=\alpha e_{31}+\beta e_{32}\\ 
C_{11}= qe_{11}+e_{22}+e_{33}+e_{44}\\
C_{22}= qe_{11}+q^2e_{22}+qe_{33}+e_{44}
\end{array}
$
&&
$\Re =\left(\matrix{*&0&0&0\cr 0&*&0&0\cr
*&*&*&*\cr0&0&0&*\cr}\right)$
        &\cr
height5pt&\omit&& \omit&&\omit&\cr
\noalign{\hrule }
height5pt&\omit&& \omit&&\omit&\cr    
&$\matrix{dim \Re \cr 7\cr
dim I \cr 1 \cr }$ && 
$I\cong C$
           && 
Perturbation zero &\cr
height5pt&\omit&& \omit&& \omit&\cr
\noalign{\hrule }
height2pt&\omit&& \omit&& \omit&\cr
\noalign{\hrule }
height5pt&\omit&& \omit&& \omit&\cr

	&{\bf CASE 3.14)}   
&&
$\begin{array}{l}
C_{12}=\gamma e_{34}\\ 
C_{21}=\beta e_{32}\\ 
C_{11}= qe_{11}+qe_{22}+e_{33}+e_{44}+e_{12}\\
C_{22}= qe_{11}+qe_{22}+qe_{33}+e_{44}-e_{12}
\end{array}
$
&&
$\Re =\left(\matrix{\epsilon&*&0&0\cr 0&\epsilon&0&0\cr
0&*&*&*\cr0&0&0&*\cr}\right)$
        &\cr
height5pt&\omit&& \omit&&\omit&\cr
\noalign{\hrule }
height5pt&\omit&& \omit&&\omit&\cr    
&$\matrix{dim \Re \cr 6\cr
dim I \cr 1 \cr }$ && 
$I\cong C$
           && 
Perturbation zero &\cr
height5pt&\omit&& \omit&& \omit&\cr
\noalign{\hrule }
height2pt&\omit&& \omit&& \omit&\cr
\noalign{\hrule }
height5pt&\omit&& \omit&& \omit&\cr
	&{\bf CASE 3.15)}   
&&
$\begin{array}{l}
C_{12}=\beta e_{23}\\ 
C_{21}=\gamma e_{43}\\ 
C_{11}= qe_{11}+qe_{22}+qe_{33}+e_{44}\\
C_{22}= qe_{11}+qe_{22}+e_{33}+e_{44}
\end{array}
$
&&
$\Re =\left(\matrix{\epsilon&0&0&0\cr 0&\epsilon&*&0\cr
0&0&*&0\cr0&0&*&*\cr}\right)$
        &\cr
height5pt&\omit&& \omit&&\omit&\cr
\noalign{\hrule }
height5pt&\omit&& \omit&&\omit&\cr    
&$\matrix{dim \Re \cr 5\cr
dim I \cr 3 \cr }$ && 
$I =\left(\matrix{\epsilon&0&0&0\cr \varphi&\alpha&0&0\cr 0&0&
\alpha& 0\cr 0&0&0&\alpha \cr}\right)$
           && 
Perturbation zero &\cr
height5pt&\omit&& \omit&& \omit&\cr
\noalign{\hrule }
height2pt&\omit&& \omit&& \omit&\cr
\noalign{\hrule }}}
\vbox{\offinterlineskip
\hrule
\halign{&\vrule#&\strut \quad \hfil# \hfil
 &\vrule#&\quad \hfil# \hfil&\vrule#& \quad# \hfil \cr
height5pt&\omit&& \omit&& \omit&\cr
	&{\bf CASE 3.16)}   
&&
$\begin{array}{l}
C_{12}=\beta e_{23} \\ 
C_{21}=\gamma e_{43}\\ 
C_{11}= \alpha e_{11}+qe_{22}+qe_{33}+e_{44}\\
C_{22}= \frac{q^2}{\alpha}e_{11}+qe_{22}+e_{33}+e_{44}\\
				\\
\alpha\neq q  
\end{array}
$
&&
$\Re =\left(\matrix{*&0&0&0\cr 
0&*&*&0\cr 0&0&*&0\cr 0&0&*&*\cr}
\right)$
        &\cr
height5pt&\omit&& \omit&&\omit&\cr
\noalign{\hrule }
height5pt&\omit&& \omit&&\omit&\cr    
&$\matrix{dim \Re \cr 6\cr 
dim I \cr 2 \cr }$ &&
$I =\left(\matrix{\alpha&0&0&0\cr 0&\beta&0&0\cr 0&0&
\beta& 0\cr 0&0&0&\beta \cr}\right)$
           && 
Perturbation zero&\cr
height5pt&\omit&& \omit&& \omit&\cr
\noalign{\hrule }
height2pt&\omit&& \omit&& \omit&\cr
\noalign{\hrule }
height5pt&\omit&& \omit&& \omit&\cr

	&{\bf CASE 3.17)}   
&&
$\begin{array}{l}
C_{12}=0\\ 
C_{21}=\gamma e_{43}\\ 
C_{11}= \alpha e_{11}+\alpha e_{22}+qe_{33}+e_{44}+e_{12}\\
C_{22}= \frac{q^2}{\alpha}e_{11}+ \frac{q^2}{\alpha}e_{22}+
e_{33}+e_{44}-\frac{q^2}{\alpha^2}e_{12}
\end{array}
$
&&
$\Re =\left(\matrix{\epsilon&*&0& 0\cr 0&\epsilon&0&0\cr
0&0&*&0\cr0&0&*&*\cr}\right)$
        &\cr
height5pt&\omit&& \omit&&\omit&\cr
\noalign{\hrule }
height5pt&\omit&& \omit&&\omit&\cr    
&$\matrix{dim \Re \cr 5\cr
dim I \cr 3 \cr }$ && 
$I =\left(\matrix{\beta&\varphi&0&0\cr 0&\beta&0&0\cr 0&0&
\alpha& 0\cr 0&0&0&\alpha \cr}
\right)$
           && 
Perturbation zero &\cr
height5pt&\omit&& \omit&& \omit&\cr
\noalign{\hrule }
height2pt&\omit&& \omit&& \omit&\cr
\noalign{\hrule }
height5pt&\omit&& \omit&& \omit&\cr

	&{\bf CASE 3.18)}   
&&
$\begin{array}{l}
C_{12}=qe_{13}\\ 
C_{21}=\gamma e_{43}\\ 
C_{11}= qe_{11}+qe_{22}+qe_{33}+e_{44}\\
C_{22}= qe_{11}+qe_{22}+e_{33}+e_{44}
\end{array}
$
&&
$\Re =\left(\matrix{\epsilon&0&*&0\cr 0&\epsilon&0&0\cr
0&0&*&0\cr0&0&0&*\cr}\right)$
        &\cr
height5pt&\omit&& \omit&&\omit&\cr
\noalign{\hrule }
height5pt&\omit&& \omit&&\omit&\cr    
&$\matrix{dim \Re \cr 4 \cr 
dim I \cr 3 \cr }$ && 
$I =\left(\matrix{\alpha&\beta&0&0\cr 0&\gamma&0&0\cr 0&0&
\alpha & 0\cr 0&0&0&\alpha \cr}\right)$
           && 
Perturbation zero &\cr
height5pt&\omit&& \omit&& \omit&\cr
\noalign{\hrule }
height2pt&\omit&& \omit&& \omit&\cr
\noalign{\hrule }}}
\vbox{\offinterlineskip
\hrule
\halign{&\vrule#&\strut \quad \hfil# \hfil
 &\vrule#&\quad \hfil# \hfil&\vrule#& \quad# \hfil \cr
\noalign{\hrule }

	&{\bf CASE 3.19)}   
&&
$\begin{array}{l}
C_{12}=\alpha e_{13}\\ 
C_{21}=\gamma e_{43}\\ 
C_{11}= qe_{11}+\alpha e_{22}+qe_{33}+e_{44}\\
C_{22}= qe_{11}+\frac{q^2}{\alpha}e_{22}+e_{33}+e_{44}\\
			\\
\alpha\neq q
\end{array}
$
&&
$\Re =\left(\matrix{*&0&*&0\cr 0&*&0&0\cr
0&0&*&0\cr 0&0&*&*\cr}\right)$
        &\cr
height5pt&\omit&& \omit&&\omit&\cr
\noalign{\hrule }
height5pt&\omit&& \omit&&\omit&\cr    
&$\matrix{dim \Re \cr 6\cr
dim I \cr 2 \cr }$ && 
$I =\left(\matrix{\alpha&0&0&0\cr 0&\beta&0&0\cr 0&0&
\alpha & 0\cr 0&0&0&\alpha \cr}\right)$
           && 
Perturbation zero&\cr
height5pt&\omit&& \omit&& \omit&\cr
\noalign{\hrule }
height2pt&\omit&& \omit&& \omit&\cr
\noalign{\hrule }
height5pt&\omit&& \omit&& \omit&\cr


	&{\bf CASE 3.20)}   
&&
$\begin{array}{l}
C_{12}=\alpha e_{13} \\ 
C_{21}=\gamma e_{43} \\ 
C_{11}= qe_{11}+qe_{22}+qe_{33}+e_{44}+e_{12}\\
C_{22}= qe_{11}+qe_{22}+e_{33}+e_{44}-e_{12}
\end{array}
$
&&
$\Re =\left(\matrix{\epsilon&*&*&0\cr 0&\epsilon&0&0\cr
0&0&*&*\cr0&0&0&*\cr}\right)$
        &\cr
height5pt&\omit&& \omit&&\omit&\cr
\noalign{\hrule }
height2pt&\omit&& \omit&&\omit&\cr    
&$\matrix{dim \Re \cr 6\cr
dim I \cr 2 \cr }$ && 
$I =\left(\matrix{\alpha&\beta&0&0\cr 0&\alpha&0&0\cr 0&0&
\alpha & 0\cr 0&0&0&\alpha \cr}\right)$
           && 
Perturbation zero &\cr
height5pt&\omit&& \omit&& \omit&\cr
\noalign{\hrule }
height2pt&\omit&& \omit&& \omit&\cr
\noalign{\hrule }}}

\centerline{\Huge  CASE 4) $d$=$diag(q^3,q^2,q,1)$}
\vspace*{1.5cm}
\vbox{\offinterlineskip
\hrule
\halign{&\vrule#&\strut \quad \hfil# \hfil
 &\vrule#&\quad \hfil# \hfil&\vrule#& \quad# \hfil \cr
height5pt&\omit&& \omit&& \omit&\cr&
        {\bf CASE 4.1)}   
&&
$\begin{array}{l}
C_{12}=0\\ 
C_{21}=\alpha e_{21}+\beta e_{32}+\gamma e_{43} \\
C_{11}=q^3 e_{11}+q^2 e_{22}+qe_{33}+e_{44}\\
C_{22}={\bf 1}\\
\end{array}$
&& 
        $\Re =\left(\matrix{*&0&0&0\cr *&*&0&0\cr 0&*&*&0
\cr 0&0&*&*\cr}
\right)$&\cr
height5pt&\omit&& \omit&&\omit&\cr
\noalign{\hrule }
height5pt&\omit&& \omit&&\omit&\cr    
&$\matrix{dim \Re \cr 7 \cr 
dim I \cr 1}$ &&
         $\begin{array}{l}
I \cong C \end{array}$ && 
         Perturbation zero &\cr
height5pt&\omit&& \omit&& \omit&\cr
\noalign{\hrule }
height2pt&\omit&& \omit&& \omit&\cr
\noalign{\hrule }
height5pt&\omit&& \omit&& \omit&\cr&
        {\bf CASE 4.2)}               
&&
$\begin{array}{c}
C_{12}=\beta e_{23}\\
C_{21}=\alpha e_{21}+\gamma e_{43}\\
C_{11}=q^2e_{11}+qe_{22}+qe_{33}+e_{44}\\
C_{22}=qe_{11}+qe_{22}+e_{33}+e_{44}
\end{array}
$
&& 
        $\Re =\left(\matrix{*&0&0&0\cr *&*&*&0\cr 0&0&
*&0\cr 0&0&*&*\cr}
\right)$&\cr
height5pt&\omit&& \omit&&\omit&\cr
\noalign{\hrule }
height5pt&\omit&& \omit&&\omit&\cr    
&
$\matrix{dim \Re \cr 7 \cr dim I \cr 1}$
 &&
$I\cong C$
&& 
         Perturbation zero &\cr
height5pt&\omit&& \omit&& \omit&\cr
\noalign{\hrule }
height2pt&\omit&& \omit&& \omit&\cr
\noalign{\hrule }
height5pt&\omit&& \omit&& \omit&\cr&
{\bf CASE 4.3)}   
&&
$\begin{array}{c}
	C_{12}=\alpha e_{12}+\gamma e_{34}\\
	C_{21}=\beta e_{32}\\
	C_{11}=qe_{11}+qe_{22}+e_{33}+e_{44}\\
        C_{22}=q^{2}e_{11}+qe_{22}+qe_{33}+e_{44}
\end{array}$
&& 
        $\Re =\left(\matrix{*&*&0&0\cr 0&*&0&0\cr 0&*&
*&*\cr 0&0&0&*\cr}
\right)$&\cr
height5pt&\omit&& \omit&&\omit&\cr
\noalign{\hrule }
height5pt&\omit&& \omit&&\omit&\cr    
&$\matrix{dim \Re \cr 7 \cr 
dim I \cr 1}$ 
         &&
$I \cong C$
&& 
         Perturbation zero&\cr
height5pt&\omit&& \omit&& \omit&\cr
\noalign{\hrule }
height2pt&\omit&& \omit&& \omit&\cr
\noalign{\hrule }}}
\vbox{\offinterlineskip
\hrule
\halign{&\vrule#&\strut \quad \hfil# \hfil
 &\vrule#&\quad \hfil# \hfil&\vrule#& \quad# \hfil \cr
height5pt&\omit&& \omit&& \omit&\cr
        &{\bf CASE 4.4)}  
&&
$\begin{array}{l}
C_{12}=\alpha e_{12}+\beta e_{23}+\gamma e_{34}\\
C_{21}=0\\
C_{11}={\bf 1}\\ 
C_{22}=q^{3}e_{11}+q^{2}e_{22}+qe_{33}+e_{44}
\end{array}
$
&& 
        $\Re =\left(\matrix{*&*&0&0\cr 0&*&*&0\cr 0&0&
*&*\cr 0&0&0&*\cr}
\right)$&\cr
height5pt&\omit&& \omit&&\omit&\cr
\noalign{\hrule }
height5pt&\omit&& \omit&&\omit&\cr    
&$\matrix{dim \Re \cr 7 \cr 
dim I \cr 1}$ && 
$I \cong C$&& 
Perturbation zero &\cr
height5pt&\omit&& \omit&& \omit&\cr
\noalign{\hrule }
height2pt&\omit&& \omit&& \omit&\cr
\noalign{\hrule }
height5pt&\omit&& \omit&& \omit&\cr
        & {\bf CASE 4.5)}   
&&
$\begin{array}{c}
	C_{12}=\gamma e_{34}\\
	C_{21}=\alpha e_{21}+\beta e_{32}\\
	C_{11}=q^2e_{11}+qe_{22}+e_{33}+e_{44}\\
        C_{22}=qe_{11}+qe_{22}+qe_{33}+e_{44}
\end{array}$
&& 
        $\Re =\left(\matrix{*&0&0&0\cr *
&*&0&0\cr 0&*&*&*\cr 0&0&0&*\cr}
\right)$&\cr
height5pt&\omit&& \omit&&\omit&\cr
\noalign{\hrule }
height5pt&\omit&& \omit&&\omit&\cr    
&$\matrix{dim \Re \cr 7\cr 
dim I \cr 1}$ && 
$I =C$
           && 
Perturbation zero &\cr
height5pt&\omit&& \omit&& \omit&\cr
\noalign{\hrule }
height2pt&\omit&& \omit&& \omit&\cr
\noalign{\hrule }
height5pt&\omit&& \omit&& \omit&\cr
        & {\bf CASE 4.6)}   
&&
$\begin{array}{l}
C_{12}=\beta e_{23}+\gamma e_{34} \\ 
C_{21}=\alpha e_{21} \\ 
C_{11}=q e_{11}+e_{22}+e_{33}+e_{44}\\
C_{22}= q^2e_{11}+q^{2}e_{22}+qe_{33}+e_{44}
\end{array}
$
&& 
        $\Re =\left(\matrix{*&0&0&0\cr *&*&*&0\cr 0&0&
* &*\cr 0&0&0&*\cr}
\right)$&\cr
height5pt&\omit&& \omit&&\omit&\cr
\noalign{\hrule }
height5pt&\omit&& \omit&&\omit&\cr    
&$\matrix{dim \Re \cr 7\cr 
dim I \cr 1}$ && 
$I\cong C$
           && 
Perturbation zero &\cr
height5pt&\omit&& \omit&& \omit&\cr
\noalign{\hrule }
height2pt&\omit&& \omit&& \omit&\cr
\noalign{\hrule }
height5pt&\omit&& \omit&& \omit&\cr
        & {\bf CASE 4.7)}   
&&
$\begin{array}{l}
C_{12}=\alpha e_{12} \\ 
C_{21}=\beta e_{32}+\gamma e_{43} \\ 
C_{11}= q^2e_{11}+q^2e_{22}+qe_{33}+e_{44}\\
C_{22}= qe_{11}+e_{22}+e_{33}+e_{44}
\end{array}
$
&& 
        $\Re =\left(\matrix{*&*&0&0\cr 0&*&0&0\cr 0&*&* &0\cr 
0&0&*&* \cr}
\right)$&\cr
height5pt&\omit&& \omit&&\omit&\cr
\noalign{\hrule }
height5pt&\omit&& \omit&&\omit&\cr    
&$\matrix{dim \Re \cr 7 \cr 
dim I \cr 1}$ && 
$I\cong C$
           && 
Perturbation zero &\cr
height5pt&\omit&& \omit&& \omit&\cr
\noalign{\hrule }
height2pt&\omit&& \omit&& \omit&\cr
\noalign{\hrule }}}
\vbox{\offinterlineskip
\hrule
\halign{&\vrule#&\strut \quad \hfil# \hfil
 &\vrule#&\quad \hfil# \hfil&\vrule#& \quad# \hfil \cr
height5pt&\omit&& \omit&& \omit&\cr
	&{\bf CASE 4.8)}   
&&
$\begin{array}{l}
C_{12}=\alpha e_{12}+\beta e_{23} \\ 
C_{21}=\gamma e_{43}\\ 
C_{11}= qe_{11}+qe_{22}+qe_{33}+e_{44}\\
C_{22}= q^2e_{11}+qe_{22}+e_{33}+e_{44}
\end{array}
$
&& 
        $\Re =\left(\matrix{*&*&0&0\cr 0&*&*&0\cr 0&0&*&0\cr 
0&0&*&* \cr}
\right)$&\cr
height5pt&\omit&& \omit&&\omit&\cr
\noalign{\hrule }
height5pt&\omit&& \omit&&\omit&\cr    
&$\matrix{dim \Re \cr 7\cr 
dim I \cr 1 \cr }$ && 
$I \cong C$
           && 
Perturbation zero &\cr
height5pt&\omit&& \omit&& \omit&\cr
\noalign{\hrule }}}
\newpage
\centerline{\Huge  CASE 5) $d$=$diag(\alpha,q^2,q,1)$}
\centerline{}
\centerline{\Large $\alpha\neq 0, q^{-1},1,q,q^2,q^3$}
\vspace*{1.5cm}
\vbox{\offinterlineskip
\hrule
\halign{&\vrule#&\strut \quad \hfil# \hfil
 &\vrule#&\quad \hfil# \hfil&\vrule#& \quad# \hfil \cr
height5pt&\omit&& \omit&& \omit&\cr&

        {\bf CASE 5.1)}   
&&
$\begin{array}{l}
C_{12}=0\\ 
C_{21}=\gamma e_{32}+\beta e_{43} \\
C_{11}=\delta e_{11}+q^2 e_{22}+qe_{33}+e_{44}\\
C_{22}=\frac{\alpha}{\delta}e_{11}+e_{22}+e_{33}+e_{44}\\
\end{array}$
&& 
        $\Re =\left(\matrix{*&0&0&0\cr 0&*&0&0\cr 0&*&*&0
\cr 0&0&*&*\cr}
\right)$&\cr
height5pt&\omit&& \omit&&\omit&\cr
\noalign{\hrule }
height5pt&\omit&& \omit&&\omit&\cr    
&$\matrix{dim \Re \cr 6 \cr 
dim I \cr 2}$ &&
$I =\left(\matrix{\beta&0&0&0\cr 0&\gamma&0&0\cr 0&0&
\gamma & 0\cr 0&0&0&\gamma \cr}
\right)$
       && 
         Perturbation zero &\cr
height5pt&\omit&& \omit&& \omit&\cr
\noalign{\hrule }
height2pt&\omit&& \omit&& \omit&\cr
\noalign{\hrule }
height5pt&\omit&& \omit&& \omit&\cr& 

        {\bf CASE 5.2)}               
&&
$\begin{array}{c}
C_{12}=\beta e_{23}+\gamma e_{34}\\
C_{21}=0\\
C_{11}=\delta e_{11}+e_{22}+e_{33}+e_{44}\\
C_{22}=\frac{\alpha}{\delta}e_{11}+q^2e_{22}+qe_{33}+e_{44}
\end{array}
$
&& 
        $\Re =\left(\matrix{*&0&0&0\cr 0&*&*&0\cr 0&0&
*&*\cr 0&0&0&* \cr}
\right)$&\cr
height5pt&\omit&& \omit&&\omit&\cr
\noalign{\hrule }
height5pt&\omit&& \omit&&\omit&\cr    
&
$\matrix{dim \Re \cr 6 \cr dim I \cr 2}$
 &&
$I =\left(\matrix{\gamma&0&0&0\cr 0&\beta&0&0\cr 0&0&
\beta & 0\cr 0&0&0&\beta \cr}
\right)$
&& 
         Perturbation zero &\cr
height5pt&\omit&& \omit&& \omit&\cr
\noalign{\hrule }
height2pt&\omit&& \omit&& \omit&\cr
\noalign{\hrule }
height5pt&\omit&& \omit&& \omit&\cr&
{\bf CASE 5.3)}   
&& 
$\begin{array}{c}
	C_{12}=\beta e_{34}\\
	C_{21}=\gamma e_{32}\\
	C_{11}=\delta e_{11}+qe_{22}+e_{33}+e_{44}\\
        C_{22}=\frac{\alpha}{\delta}e_{11}+qe_{22}+qe_{33}+e_{44}
\end{array}$
&& 
        $\Re =\left(\matrix{*&0&0&0\cr 0&*&0&0\cr 0&*&
*&*\cr 0&0&0&*\cr}
\right)$&\cr
height5pt&\omit&& \omit&&\omit&\cr
\noalign{\hrule }
height5pt&\omit&& \omit&&\omit&\cr    
&$\matrix{dim \Re \cr 6 \cr 
dim I \cr 2}$ 
         &&
$I =\left(\matrix{\gamma&0&0&0\cr 0&\beta&0&0\cr 0&0&
\beta & 0\cr 0&0&0&\beta \cr}
\right)$
&& 
         Perturbation zero&\cr
height5pt&\omit&& \omit&& \omit&\cr
\noalign{\hrule }
height2pt&\omit&& \omit&& \omit&\cr
\noalign{\hrule }}}
\vbox{\offinterlineskip
\hrule
\halign{&\vrule#&\strut \quad \hfil# \hfil
 &\vrule#&\quad \hfil# \hfil&\vrule#& \quad# \hfil \cr
height5pt&\omit&& \omit&& \omit&\cr

        &{\bf CASE 5.4)}  
&&
$\begin{array}{l}
C_{12}=\gamma e_{23}\\
C_{21}=\beta e_{43}\\
C_{11}=\delta e_{11}+qe_{22}+qe_{33}+e_{44}\\ 
C_{22}=\frac{\alpha}{\delta}e_{11}+qe_{22}+e_{33}+e_{44}
\end{array}
$
&& 
        $\Re =\left(\matrix{*&0&0&0\cr 0&*&*&0\cr 0&0&
*&0\cr 0&0&*&*\cr}
\right)$&\cr
height5pt&\omit&& \omit&&\omit&\cr
\noalign{\hrule }
height5pt&\omit&& \omit&&\omit&\cr    
&$\matrix{dim \Re \cr 6 \cr 
dim I \cr 2}$ && 
$I =\left(\matrix{\gamma&0&0&0\cr 0&\beta&0&0\cr 0&0&
\beta & 0\cr 0&0&0&\beta \cr}
\right)$
&& 
Perturbation zero &\cr
height5pt&\omit&& \omit&& \omit&\cr
\noalign{\hrule }
height2pt&\omit&& \omit&& \omit&\cr
\noalign{\hrule }}}
\

\centerline{\Huge  CASE 6) $d$=$diag(q^2,q^2,q,1)+e_{12}$}
\vspace*{1.5cm}
\vbox{\offinterlineskip
\hrule
\halign{&\vrule#&\strut \quad \hfil# \hfil
 &\vrule#&\quad \hfil# \hfil&\vrule#& \quad# \hfil \cr
height5pt&\omit&& \omit&& \omit&\cr&
        {\bf CASE 6.1)}   
&&
$\begin{array}{l}
C_{12}=0\\ 
C_{21}=\alpha e_{43}\\
C_{11}=q^2\beta e_{11}+q^2\beta e_{22}+q\gamma e_{33}+\gamma e_{44}+
\beta e_{12} \\
c_{22}=\beta^{-1}e_{11}+\beta^{-1}e_{22}+q\delta^{-1}e_{33}+
\gamma^{-1}e_{44}\\
\end{array}$
&& 
        $\Re =\left(\matrix{\epsilon&*&0&0\cr 0&\epsilon&0&0\cr 0&0&
*&0\cr 0&0&*&*\cr}
\right)$&\cr
height5pt&\omit&& \omit&&\omit&\cr
\noalign{\hrule }
height5pt&\omit&& \omit&&\omit&\cr    
&$\matrix{dim \Re \cr 5 \cr 
dim I \cr 3}$ &&
$I =\left(\matrix{\beta&\gamma&0&0\cr 0&\beta&0&0\cr 0&0&
\alpha & 0\cr 0&0&0&\alpha \cr}
\right)$         
&& 
         Perturbation zero &\cr
height5pt&\omit&& \omit&& \omit&\cr
\noalign{\hrule }
height2pt&\omit&& \omit&& \omit&\cr
\noalign{\hrule }
height5pt&\omit&& \omit&& \omit&\cr&
        {\bf CASE 6.2)}               
&&
$\begin{array}{c}
C_{12}=\alpha e_{13}+\beta e_{34}\\
C_{21}=0\\
C_{11}=q^2e_{11}+q^2e_{22}+q^2e_{33}+q^2e_{44}+e_{12}\\
C_{22}=e_{11}+e_{22}+q^{-1}e_{33}+q^{-2}e_{44}
\end{array}
$
&& 
        $\Re =\left(\matrix{\epsilon&*&*&0\cr 0&\epsilon&0&0\cr 0&0&
* &*\cr 0&0&0&* \cr}
\right)$&\cr
height5pt&\omit&& \omit&&\omit&\cr
\noalign{\hrule }
height5pt&\omit&& \omit&&\omit&\cr    
&
$\matrix{dim \Re \cr 6 \cr dim I \cr 2}$
 &&
$I =\left(\matrix{\alpha&\beta&0&0\cr 0&\alpha&0&0\cr 0&0&
\alpha & 0\cr 0&0&0&\alpha \cr}
\right)$
&& 
         Perturbation zero &\cr
height5pt&\omit&& \omit&& \omit&\cr
\noalign{\hrule }
height2pt&\omit&& \omit&& \omit&\cr
\noalign{\hrule }
height5pt&\omit&& \omit&& \omit&\cr&
{\bf CASE 6.3)}   
&&
$\begin{array}{c}
	C_{12}=\alpha e_{34}\\
	C_{21}=\beta e_{32}\\
	C_{11}=q^2e_{11}+q^{2}e_{22}+qe_{33}+qe_{44}+e_{12}\\
        C_{22}=e_{11}+e_{22}+e_{33}+q^{-1}e_{44}
\end{array}$
&& 
        $\Re =\left(\matrix{\epsilon&*&0&0\cr 0&\epsilon&0&0\cr 0&*&
*&0\cr 0&0&*&* \cr}
\right)$&\cr
height5pt&\omit&& \omit&&\omit&\cr
\noalign{\hrule }
height5pt&\omit&& \omit&&\omit&\cr    
&$\matrix{dim \Re \cr 6 \cr 
dim I \cr 2}$ 
         &&
$I =\left(\matrix{\alpha&\beta&0&0\cr 0&\alpha&0&0\cr 0&0&
\alpha & 0\cr 0&0&0&\alpha \cr}
\right)$
&& 
         Perturbation zero&\cr
height5pt&\omit&& \omit&& \omit&\cr
\noalign{\hrule }
height2pt&\omit&& \omit&& \omit&\cr
\noalign{\hrule }}}
\vbox{\offinterlineskip
\hrule
\halign{&\vrule#&\strut \quad \hfil# \hfil
 &\vrule#&\quad \hfil# \hfil&\vrule#& \quad# \hfil \cr
height5pt&\omit&& \omit&& \omit&\cr
        &{\bf CASE 6.4)}  
&&
$\begin{array}{l}
C_{12}=\alpha e_{34}\\
C_{21}=0\\
C_{11}=q^2\beta e_{11}+q^2\beta e_{22}+\gamma e_{33}+\gamma e_{44}
+\beta e_{12}\\ 
C_{22}=\beta^{-1}e_{11}+\beta^{-1}e_{22}+q\gamma^{-1}e_{33}+
\gamma^{-1}e_{44}
\end{array}
$
&& 
        $\Re =\left(\matrix{\epsilon&*&0&0\cr 0&\epsilon&0&0\cr 0&0&
*&*\cr 0&0&0&*\cr}
\right)$&\cr
height5pt&\omit&& \omit&&\omit&\cr
\noalign{\hrule }
height5pt&\omit&& \omit&&\omit&\cr    
&$\matrix{dim \Re \cr 5 \cr 
dim I \cr 3}$ && 
$I =\left(\matrix{\alpha&\gamma&0&0\cr 0&\alpha&0&0\cr 0&0&
\beta & 0\cr 0&0&0&\beta \cr}
\right)$&& 
Perturbation zero &\cr
height5pt&\omit&& \omit&& \omit&\cr
\noalign{\hrule }
height2pt&\omit&& \omit&& \omit&\cr
\noalign{\hrule }
height5pt&\omit&& \omit&& \omit&\cr
        & {\bf CASE 6.5)}   
&&
$\begin{array}{c}
	C_{12}=\alpha e_{13}\\
	C_{21}=0\\
	C_{11}=q^2\beta e_{11}+q^2\beta e_{22}+q^2\beta e_{33}
+\delta e_{44}+\beta e_{12}\\
        C_{22}=\beta^{-1}e_{11}+\beta^{-1}e_{22}+q^{-1}\beta^{-1}e_{33}+
\delta^{-1}e_{44}
\end{array}$
&& 
        $\Re =\left(\matrix{\epsilon&*&*&0\cr 0
&\epsilon&0&0\cr 0&0&
*&0\cr 0&0&0&*\cr}
\right)$&\cr
height5pt&\omit&& \omit&&\omit&\cr
\noalign{\hrule }
height5pt&\omit&& \omit&&\omit&\cr    
&$\matrix{dim \Re \cr 5 \cr 
dim I \cr 3}$ && 
$I =\left(\matrix{\alpha&\gamma&0&0\cr 0&\alpha&0&0\cr 0&0&
\alpha& 0\cr 0&0&0&\beta \cr}
\right)$
           && 
Perturbation zero &\cr
height5pt&\omit&& \omit&& \omit&\cr
\noalign{\hrule }
height2pt&\omit&& \omit&& \omit&\cr
\noalign{\hrule }
height5pt&\omit&& \omit&& \omit&\cr
        & {\bf CASE 6.6)}   
&&
$\begin{array}{l}
C_{12}=0 \\ 
C_{21}=\alpha e_{32} \\ 
C_{11}= q^2\beta e_{11}+q^2\beta e_{22}+q\beta e_{33}+
\delta e_{44}+\beta e_{12}\\
C_{22}= \beta^{-1}e_{11}+\beta^{-1}e_{22}+\beta^{-1}e_{33}
+\delta^{-1}e_{44}
\end{array}
$
&& 
        $\Re =\left(\matrix{\epsilon&*&0&0\cr 0&\epsilon&0&0\cr 0&*&
*&0\cr 0&0&0&*\cr}
\right)$&\cr
height5pt&\omit&& \omit&&\omit&\cr
\noalign{\hrule }
height5pt&\omit&& \omit&&\omit&\cr    
&$\matrix{dim \Re \cr 5 \cr 
dim I \cr 3}$ && 
$I =\left(\matrix{\alpha&\gamma&0&0\cr 0&\alpha&0&0\cr 0&0&
\alpha & 0\cr 0&0&0&\beta \cr}
\right)$
           && 
Perturbation zero &\cr
height5pt&\omit&& \omit&& \omit&\cr
\noalign{\hrule }
height2pt&\omit&& \omit&& \omit&\cr
\noalign{\hrule }
height5pt&\omit&& \omit&& \omit&\cr
        & {\bf CASE 6.7)}   
&&
$\begin{array}{l}
C_{12}=0 \\ 
C_{21}=0\\ 
C_{11}= q^2\alpha e_{11}+q^2\alpha e_{22}+\beta e_{33}+\gamma e_{44}
+\alpha e_{12}\\
C_{22}= \alpha^{-1}e_{11}+\alpha^{-1}e_{22}+q\beta^{-1}e_{33}+
\gamma^{-1}e_{44}
\end{array}
$
&& 
        $\Re =\left(\matrix{\epsilon&*&0&0\cr 0&\epsilon&0&0\cr 0&0&
*&0\cr 0&0&0&*\cr}
\right)$&\cr
height5pt&\omit&& \omit&&\omit&\cr
\noalign{\hrule }
height5pt&\omit&& \omit&&\omit&\cr    
&$\matrix{dim \Re \cr 4 \cr 
dim I \cr 4}$ && 
$I =\left(\matrix{\gamma&\varphi&0&0\cr 0&\gamma&0&0\cr 0&0&
\beta & 0\cr 0&0&0&\alpha \cr}
\right)$
           && 
Perturbation zero &\cr
height5pt&\omit&& \omit&& \omit&\cr
\noalign{\hrule }
height2pt&\omit&& \omit&& \omit&\cr
\noalign{\hrule }}}

\

\centerline{\Huge  CASE 7) $d$=$diag(q^2,q,1,1)+e_{34}$}
\vspace*{1.5cm}
\vbox{\offinterlineskip
\hrule
\halign{&\vrule#&\strut \quad \hfil# \hfil
 &\vrule#&\quad \hfil# \hfil&\vrule#& \quad# \hfil \cr
height5pt&\omit&& \omit&& \omit&\cr&
        {\bf CASE 7.1)}   
&&
$\begin{array}{l}
C_{12}=\alpha e_{24}\\
C_{21}=\beta e_{21}\\ 
C_{11}=qe_{11}+e_{22}+e_{33}+e_{44}+e_{34}\\
C_{22}=qe_{11}+qe_{22}+e_{33}+e_{44}\\
\end{array}$
&& 
        $\Re =\left(\matrix{*&0&0&0\cr 
*&*&0&*\cr 
0&0&\epsilon&*\cr 
0&0&0&\epsilon \cr}
\right)$&\cr
height5pt&\omit&& \omit&&\omit&\cr
\noalign{\hrule }
height5pt&\omit&& \omit&&\omit&\cr    
&$\matrix{dim \Re \cr 6 \cr 
dim I \cr 1}$ &&
         $I\cong C$
&& 
         Perturbation zero &\cr
height5pt&\omit&& \omit&& \omit&\cr
\noalign{\hrule }
height2pt&\omit&& \omit&& \omit&\cr
\noalign{\hrule }
height5pt&\omit&& \omit&& \omit&\cr&

        {\bf CASE 7.2)}               
&&
$\begin{array}{c}
C_{12}=\alpha e_{24}\\
C_{21}=\beta e_{21}\\
C_{11}=q\beta e_{11}+\beta e_{22}+\gamma e_{33}+\beta e_{44}+\beta e_{34}\\
C_{22}=q\delta^{-1}e_{11}+q\delta^{-1}e_{22}+\gamma^{-1}e_{33}+
\delta^{-1}e_{44}\\
		\\
\beta\neq \gamma\neq \delta
\end{array}
$
&& 
        $\Re =\left(\matrix{*&0&0&0\cr 
*&*&0&*\cr 0&0&*&*\cr 0&0&0&*\cr}\right)$&\cr
height5pt&\omit&& \omit&&\omit&\cr
\noalign{\hrule }
height5pt&\omit&& \omit&&\omit&\cr    
&
$\matrix{dim \Re \cr 7 \cr dim I \cr 1}$
 &&
$I\cong C$
&& Perturbation zero&\cr
height5pt&\omit&& \omit&& \omit&\cr
\noalign{\hrule }
height2pt&\omit&& \omit&& \omit&\cr
\noalign{\hrule }
height5pt&\omit&& \omit&& \omit&\cr&

{\bf CASE 7.3)}   
&&
$\begin{array}{c}
	C_{12}=0\\
	C_{21}=\alpha e_{21}+\beta e_{32}\\
	C_{11}=q^{2}e_{11}+qe_{22}+e_{33}+e_{44}+e_{34}\\
        C_{22}={\bf 1}
\end{array}$
&& 
        $\Re =\left(\matrix{*&0&0&0\cr *&*&0&0\cr 
0&*&\epsilon&*\cr 
0&0&0&\epsilon\cr}
\right)$&\cr
height5pt&\omit&& \omit&&\omit&\cr
\noalign{\hrule }
height5pt&\omit&& \omit&&\omit&\cr    
&$\matrix{dim \Re \cr 6 \cr 
dim I \cr 1}$ 
         &&
$I\cong C$
&& 
         Perturbation zero&\cr
height5pt&\omit&& \omit&& \omit&\cr
\noalign{\hrule }
height2pt&\omit&& \omit&& \omit&\cr
\noalign{\hrule }}}
\vbox{\offinterlineskip
\hrule
\halign{&\vrule#&\strut \quad \hfil# \hfil
&\vrule#&\quad \hfil# \hfil&\vrule#& \quad# \hfil \cr
height5pt&\omit&& \omit&& \omit&\cr
&{\bf CASE 7.4)}  
&&
$\begin{array}{l}
C_{12}=0\\
C_{21}=\alpha e_{21}+\beta e_{32}\\
C_{11}=q^2\gamma e_{11}+q\gamma e_{22}+\gamma e_{33}+\delta e_{44}
+\delta e_{34}\\ 
C_{22}=\gamma^{-1}e_{11}+\gamma^{-1}e_{22}+\gamma^{-1}e_{33}+
\delta^{-1}e_{44}\\
		\\
\gamma\neq \delta
\end{array}
$
&& 
        $\Re =\left(\matrix{*&0&0&0\cr *&*&0&0\cr 0&*&
* &*\cr 0&0&0&* \cr}
\right)$&\cr
height5pt&\omit&& \omit&&\omit&\cr
\noalign{\hrule }
height5pt&\omit&& \omit&&\omit&\cr    
&$\matrix{dim \Re \cr 7 \cr 
dim I \cr 1}$ && 
$I\cong C$&& 
Perturbation zero &\cr
height5pt&\omit&& \omit&& \omit&\cr
\noalign{\hrule }
height2pt&\omit&& \omit&& \omit&\cr
\noalign{\hrule }
height5pt&\omit&& \omit&& \omit&\cr

        & {\bf CASE 7.5)}   
&&
$\begin{array}{c}
	C_{12}=0\\
	C_{21}=\alpha e_{21}\\
	C_{11}=q\beta e_{11}+\beta e_{22}+\gamma e_{33}+\delta e_{44}
+\delta e_{34}\\
        C_{22}=q\beta^{-1}e_{11}+q\beta^{-1}e_{22}+\gamma^{-1}e_{33}+
\delta^{-1}e_{44}\\
	\\
\mbox{ either }\gamma=\delta=\beta\;\mbox{ or }\;
\gamma=\delta=q\beta\\
\mbox{ or }\gamma=\delta\neq \beta\;\mbox{ or }\;
\gamma=\delta\neq q\beta  
\end{array}$
&& 
        $\Re =\left(\matrix{*&0&0&0\cr *
&*&0&0\cr 0&0&
\epsilon&*\cr 0&0&0&\epsilon\cr}
\right)$&\cr
height5pt&\omit&& \omit&&\omit&\cr
\noalign{\hrule }
height5pt&\omit&& \omit&&\omit&\cr    
&$\matrix{dim \Re \cr 5 \cr 
dim I \cr 3}$ && 
$I =\left(\matrix{\alpha&0&0&0\cr 0&\alpha&0&0\cr 
0&0&\beta&\gamma\cr 0&0&0&\beta\cr}
\right)$
           && 
Perturbation zero &\cr
height5pt&\omit&& \omit&& \omit&\cr
\noalign{\hrule }
height2pt&\omit&& \omit&& \omit&\cr
\noalign{\hrule }
height5pt&\omit&& \omit&& \omit&\cr

        & {\bf CASE 7.6)}   
&&
$\begin{array}{l}
C_{12}=0\\ 
C_{21}=\alpha e_{21}\\ 
C_{11}= q\beta e_{11}+\beta e_{22}+\gamma e_{33}+\delta e_{44}+\delta e_{34}\\
C_{22}= q\beta^{-1}e_{11}+q\beta^{-1}e_{22}+\gamma^{-1}e_{33}+
\delta^{-1}e_{44}\\
		\\
\mbox{ neither }\;\; \gamma=\delta=\beta\;\;\mbox{ nor }
\;\;\gamma=\delta=q\beta\\
\hspace*{0.5cm}\mbox{ nor }\gamma=\delta\neq \beta
\;\;\mbox{ nor }\;\;\gamma=\delta\neq q\beta
\end{array}
$
	&& 
        $\Re =\left(\matrix{*&0&0&0\cr *&*&0&0\cr 0&0&
*&*\cr 0&0&0&* \cr}
\right)$&\cr
height5pt&\omit&& \omit&&\omit&\cr
\noalign{\hrule }
height5pt&\omit&& \omit&&\omit&\cr    
&$\matrix{dim \Re \cr 6 \cr 
dim I \cr 2}$ && 
$I =\left(\matrix{\beta&0&0&0\cr 0&\beta&0&0\cr 0&0&
\alpha& 0\cr 0&0&0&\alpha \cr}
\right)$
           && 
Perturbation zero &\cr
height5pt&\omit&& \omit&& \omit&\cr
\noalign{\hrule }
height2pt&\omit&& \omit&& \omit&\cr
\noalign{\hrule }}}
\vbox{\offinterlineskip
\hrule
\halign{&\vrule#&\strut \quad \hfil# \hfil
 &\vrule#&\quad \hfil# \hfil&\vrule#& \quad# \hfil \cr
height5pt&\omit&& \omit&& \omit&\cr
        & {\bf CASE 7.7)}   
&&
$\begin{array}{l}
C_{12}=\alpha e_{24}+\beta e_{12}\\ 
C_{21}=0 \\ 
C_{11}= {\bf 1}+e_{34}\\
C_{22}= q^2e_{11}+qe_{22}+e_{33}+e_{44}
\end{array}
$
&& 
$\Re =\left(\matrix{*&*&0&0\cr 0&*&0&*\cr 0&0&
\epsilon&*\cr 0&0&0&\epsilon\cr}
\right)$
        &\cr
height5pt&\omit&& \omit&&\omit&\cr
\noalign{\hrule }
height5pt&\omit&& \omit&&\omit&\cr    
&$\matrix{dim \Re \cr 6\cr 
dim I \cr 2}$ && 
$I =\left(\matrix{\alpha&0&0&0\cr 0&\alpha&0&0\cr 0&0&
\alpha&\beta\cr 0&0&0&\alpha \cr}
\right)$
           && 
Perturbation zero &\cr
height5pt&\omit&& \omit&& \omit&\cr
\noalign{\hrule }
height2pt&\omit&& \omit&& \omit&\cr
\noalign{\hrule }
height5pt&\omit&& \omit&& \omit&\cr

	&{\bf CASE 7.8)}   
&&
$\begin{array}{l}
C_{12}=\alpha e_{24}+\beta e_{12}\\ 
C_{21}=0\\ 
C_{11}= \alpha e_{11}+\alpha e_{22}+\gamma e_{33}+\alpha e_{44}
+\alpha e_{34}\\
C_{22}= \frac{q^2}{\alpha}e_{11}+\frac{q}{\alpha}e_{22}+\gamma^{-1}e_{33}
+\alpha^{-1}e_{44}\\
\alpha\neq \gamma
\end{array}
$
&&
$\Re =\left(\matrix{*&*&0&0\cr 
0&*&0&*\cr 0&0&*&0\cr 0&0&0&*\cr}
\right)$
        &\cr
height5pt&\omit&& \omit&&\omit&\cr
\noalign{\hrule }
height5pt&\omit&& \omit&&\omit&\cr    
&$\matrix{dim \Re \cr 6\cr 
dim I \cr 1 \cr }$ && 
$I \cong C$
           && 
Perturbation zero&\cr
height5pt&\omit&& \omit&& \omit&\cr
\noalign{\hrule }
height2pt&\omit&& \omit&& \omit&\cr
\noalign{\hrule }
height5pt&\omit&& \omit&& \omit&\cr

	&{\bf CASE 7.9)}   
&&
$\begin{array}{l}
C_{12}=\alpha e_{12}\\ 
C_{21}=\beta e_{32}\\ 
C_{11}= qe_{11}+qe_{22}+e_{33}+e_{44}+e_{34}\\
C_{22}= qe_{11}+e_{22}+e_{33}+e_{44}\end{array}
$
&&
$\Re =\left(\matrix{*&*&0& 0\cr 0&*&0&0\cr
0&*&\epsilon&*\cr0&0&0&\epsilon\cr}\right)$
        &\cr
height5pt&\omit&& \omit&&\omit&\cr
\noalign{\hrule }
height5pt&\omit&& \omit&&\omit&\cr    
&$\matrix{dim \Re \cr6\cr 
dim I \cr 1 \cr }$ && 
$I \cong C$
           && 
Perturbation zero &\cr
height5pt&\omit&& \omit&& \omit&\cr
\noalign{\hrule }
height2pt&\omit&& \omit&& \omit&\cr
\noalign{\hrule }
height5pt&\omit&& \omit&& \omit&\cr

	&{\bf CASE 7.10)}   
&&
$\begin{array}{l}
C_{12}=\alpha e_{12}\\ 
C_{21}=\beta e_{32}\\ 
C_{11}= q\gamma e_{11}+q\gamma e_{22}+\gamma e_{33}+\delta e_{44}
+\delta e_{34}\\
C_{22}= q\gamma^{-1}e_{11}+\gamma^{-1}e_{22}+\gamma^{-1}e_{33}+
\delta^{-1}e_{44}\\
\gamma\neq \delta
\end{array}
$
&&
$\Re =\left(\matrix{*&*&0& 0\cr 0&*&0&0\cr
0&*&*&*\cr0&0&0&*\cr}\right)$
        &\cr
height5pt&\omit&& \omit&&\omit&\cr
\noalign{\hrule }
height5pt&\omit&& \omit&&\omit&\cr    
&$\matrix{dim \Re \cr 7\cr 
dim I \cr 1 \cr }$ && 
$I\cong C$
           && 
Perturbation zero &\cr
height5pt&\omit&& \omit&& \omit&\cr
height5pt&\omit&& \omit&& \omit&\cr
\noalign{\hrule }
height2pt&\omit&& \omit&& \omit&\cr
\noalign{\hrule }}}
\vbox{\offinterlineskip
\hrule
\halign{&\vrule#&\strut \quad \hfil# \hfil
 &\vrule#&\quad \hfil# \hfil&\vrule#& \quad# \hfil \cr
height5pt&\omit&& \omit&& \omit&\cr

	&{\bf CASE 7.11)}   
&&
$\begin{array}{l}
C_{12}=\beta e_{12}\\ 
C_{21}=0\\ 
C_{11}= \alpha e_{11}+\alpha e_{22}+\gamma e_{33}+\gamma e_{44}+
\gamma e_{34}\\
C_{22}= \frac{q^2}{\alpha}e_{11}+\frac{q}{\alpha}e_{22}+
\gamma^{-1}e_{33}+\gamma^{-1}e_{44}
\end{array}
$
&&
$\Re =\left(\matrix{*&*&0&0\cr 0&*&0&0\cr 0&0&\epsilon&*\cr
0&0&0&\epsilon}\right)$
        &\cr
height5pt&\omit&& \omit&&\omit&\cr
\noalign{\hrule }
height5pt&\omit&& \omit&&\omit&\cr    
&$\matrix{dim \Re \cr 5\cr
dim I \cr 3 \cr }$ && 
$I =\left(\matrix{\alpha&0&0&0\cr 0&\alpha&0&0\cr 0&0&
\beta&\gamma\cr 0&0&0&\beta \cr}\right)$
          && 
Perturbation zero&\cr
height5pt&\omit&& \omit&& \omit&\cr
\noalign{\hrule }
height2pt&\omit&& \omit&& \omit&\cr
\noalign{\hrule }
height5pt&\omit&& \omit&& \omit&\cr


	&{\bf CASE 7.12)}   
&&
$\begin{array}{l}
C_{12}=\beta e_{12}\\ 
C_{21}=0 \\ 
C_{11}= \alpha e_{11}+\alpha e_{22}+\gamma e_{33}+\delta e_{44}+
\delta e_{34}\\
C_{22}= \frac{q^2}{\alpha}e_{11}+\frac{q}{\alpha}e_{22}+\gamma^{-1}e_{33}+
\delta^{-1}e_{44}\\
	\\
\gamma\neq \delta 
\end{array}
$
&&
$\Re =\left(\matrix{*&*&0&0\cr 0&*&0&0\cr
0&0&*&*\cr0&0&0&*\cr}\right)$
        &\cr
height5pt&\omit&& \omit&&\omit&\cr
\noalign{\hrule }
height2pt&\omit&& \omit&&\omit&\cr    
&$\matrix{dim \Re \cr 6\cr
dim I \cr 2 \cr }$ && 
$I =\left(\matrix{\alpha&0&0&0\cr 0&\alpha&0&0\cr 0&0&
\beta& 0\cr 0&0&0&\beta\cr}\right)$
           && 
Perturbation zero &\cr
height5pt&\omit&& \omit&& \omit&\cr
\noalign{\hrule }
height2pt&\omit&& \omit&& \omit&\cr
\noalign{\hrule }
height5pt&\omit&& \omit&& \omit&\cr

	&{\bf CASE 7.13)}   
&&
$\begin{array}{l}
C_{12}=\delta e_{24}\\ 
C_{21}=0\\ 
C_{11}= \alpha e_{11}+\beta e_{22}+\gamma e_{33}+
\beta e_{44}+\beta e_{34}\\
C_{22}= \frac{q^2}{\alpha}e_{11}+\frac{q}{\beta}e_{22}+
\gamma^{-1}e_{33}+\beta^{-1}e_{44}\\
	\\
\mbox{ either }\alpha=\beta=\gamma\;\mbox { or }\;\gamma^{-1}=q^2\alpha^{-1}=\beta^{-1}
\end{array}
$
&&
$\Re =\left(\matrix{*&0&0&0\cr 0&*&0&*\cr
0&0&\epsilon&*\cr0&0&0&\epsilon\cr}\right)$
        &\cr
height5pt&\omit&& \omit&&\omit&\cr
\noalign{\hrule }
height5pt&\omit&& \omit&&\omit&\cr    
&$\matrix{dim \Re \cr 5\cr
dim I \cr 2 \cr }$ && 
$I =\left(\matrix{\alpha&0&0&0\cr 0&\beta&0&0\cr 0&0&
\beta& 0\cr 0&0&0&\beta \cr}\right)$

           && 
Perturbation zero &\cr
height5pt&\omit&& \omit&& \omit&\cr
\noalign{\hrule }
height2pt&\omit&& \omit&& \omit&\cr
\noalign{\hrule }}}
\vbox{\offinterlineskip
\hrule
\halign{&\vrule#&\strut \quad \hfil# \hfil
 &\vrule#&\quad \hfil# \hfil&\vrule#& \quad# \hfil \cr
height5pt&\omit&& \omit&& \omit&\cr

	&{\bf CASE 7.14)}   
&&
$\begin{array}{l}
C_{12}=\delta e_{24}\\ 
C_{21}=0\\ 
C_{11}= \alpha e_{11}+\beta e_{22}+\gamma e_{33}+
\beta e_{44}+\beta e_{34}\\
C_{22}= \frac{q^2}{\alpha}e_{11}+\frac{q}{\beta}e_{22}+
\gamma^{-1}e_{33}+\beta^{-1}e_{44}\\
	\\
\mbox{ neither }\;\alpha=\beta=\gamma\;\mbox{ nor }\;
\gamma^{-1}=q^2\alpha^{-1}=\beta^{-1}
\end{array}
$
&&
$\Re =\left(\matrix{*&0&0&0\cr 0&*&0&*\cr
0&0&*&*\cr0&0&0&*\cr}\right)$
        &\cr
height5pt&\omit&& \omit&&\omit&\cr
\noalign{\hrule }
height5pt&\omit&& \omit&&\omit&\cr    
&$\matrix{dim \Re \cr 6\cr
dim I \cr 2 \cr }$ && 
$I =\left(\matrix{\alpha&0&0&0\cr 0&\beta&0&0\cr 0&0&
\beta& 0\cr 0&0&0&\beta \cr}\right)$
           && 
Perturbation zero &\cr
height5pt&\omit&& \omit&& \omit&\cr
\noalign{\hrule }
height2pt&\omit&& \omit&& \omit&\cr
\noalign{\hrule }
height5pt&\omit&& \omit&& \omit&\cr

	&{\bf CASE 7.15)}   
&&
$\begin{array}{l}
C_{12}=0\\ 
C_{21}=\beta e_{32}\\ 
C_{11}= \alpha e_{11}+q\gamma e_{22}+\gamma e_{33}+
\delta e_{44}+\delta e_{34}\\
C_{22}= \frac{q^2}{\alpha}e_{11}+\gamma^{-1}e_{22}+
\gamma^{-1}e_{33}+\delta^{-1}e_{44}\\
	\\
\mbox{ either }\alpha=\gamma=\delta\;\mbox{ or }\;\alpha\neq \gamma\neq \delta
\end{array}
$
&&
$\Re =\left(\matrix{*&0&0&0\cr 0&*&0&0\cr
0&*&\epsilon&*\cr0&0&0&\epsilon\cr}\right)$
        &\cr
height5pt&\omit&& \omit&&\omit&\cr
\noalign{\hrule }
height5pt&\omit&& \omit&&\omit&\cr    
&$\matrix{dim \Re \cr 5\cr
dim I \cr 3 \cr }$ && 
$I =\left(\matrix{\alpha&0&0&0\cr 0&\beta&0&0\cr 0&0&
\beta&\gamma\cr 0&0&0&\beta\cr}\right)$
           && 
Perturbation zero &\cr
height5pt&\omit&& \omit&& \omit&\cr
\noalign{\hrule }
height2pt&\omit&& \omit&& \omit&\cr
\noalign{\hrule }
height5pt&\omit&& \omit&& \omit&\cr
	&{\bf CASE 7.16)}   
&&
$\begin{array}{l}
C_{12}=0\\ 
C_{21}=\beta e_{32}\\ 
C_{11}= \alpha e_{11}+q\gamma e_{22}+\gamma e_{33}+
\delta e_{44}+\delta e_{34}\\
C_{22}= \frac{q^2}{\alpha}e_{11}+\gamma^{-1}e_{22}+
\gamma^{-1}e_{33}+\delta^{-1}e_{44}\\
	\\
\mbox{neither }\;\alpha=\gamma=\delta\;\mbox{ nor }\;
\alpha\neq \gamma\neq \delta

\end{array}
$
&&
$\Re =\left(\matrix{*&0&0&0\cr 
0&*&0&0\cr 0&*&*&*\cr 0&0&0&*\cr}
\right)$
        &\cr
height5pt&\omit&& \omit&&\omit&\cr
\noalign{\hrule }
height5pt&\omit&& \omit&&\omit&\cr    
&$\matrix{dim \Re \cr 6\cr 
dim I \cr 2 \cr }$ &&
$I =\left(\matrix{\alpha&0&0&0\cr 0&\beta&0&0\cr 0&0&
\beta& 0\cr 0&0&0&\beta \cr}\right)$
           && 
Perturbation zero&\cr
height5pt&\omit&& \omit&& \omit&\cr
\noalign{\hrule }
height2pt&\omit&& \omit&& \omit&\cr
\noalign{\hrule }}}
\vbox{\offinterlineskip
\hrule
\halign{&\vrule#&\strut \quad \hfil# \hfil
 &\vrule#&\quad \hfil# \hfil&\vrule#& \quad# \hfil \cr
height5pt&\omit&& \omit&& \omit&\cr

	&{\bf CASE 7.17)}   
&&
$\begin{array}{l}
C_{12}=0\\ 
C_{21}=0\\ 
C_{11}= {\bf 1}+e_{34}\\
C_{22}= q^2e_{11}+ qe_{22}+e_{33}+e_{44}
\end{array}
$
&&
$\Re =\left(\matrix{*&0&0& 0\cr 0&*&0&0\cr
0&0&\epsilon&*\cr0&0&0&\epsilon\cr}\right)$
        &\cr
height5pt&\omit&& \omit&&\omit&\cr
\noalign{\hrule }
height5pt&\omit&& \omit&&\omit&\cr    
&$\matrix{dim \Re \cr4\cr 
dim I \cr 4 \cr }$ && 
$I =\left(\matrix{\alpha&0&0&0\cr 0&\beta&0&0\cr 0&0&
\gamma& \delta\cr 0&0&0&\gamma \cr}
\right)$
           && 
Perturbation zero &\cr
height5pt&\omit&& \omit&& \omit&\cr
\noalign{\hrule }
height2pt&\omit&& \omit&& \omit&\cr
\noalign{\hrule }
height5pt&\omit&& \omit&& \omit&\cr

	&{\bf CASE 7.18)}   
&&
$\begin{array}{l}
C_{12}=0\\ 
C_{21}=0\\ 
C_{11}= \alpha e_{11}+\beta e_{22}+\gamma e_{33}+
\delta e_{44}+\delta e_{34}\\
C_{22}= \frac{q^2}{\alpha}e_{11}+\frac{q}{\beta}e_{22}+
\gamma^{-1}e_{33}+\epsilon^{-1}e_{44}
    			\\
\alpha\neq \beta\neq \gamma \neq \delta \neq \epsilon
\end{array}
$
&&
$\Re =\left(\matrix{*&0&0&0\cr 0&*&0&0\cr
0&0&*&*\cr0&0&0&*\cr}\right)$
        &\cr
height5pt&\omit&& \omit&&\omit&\cr
\noalign{\hrule }
height5pt&\omit&& \omit&&\omit&\cr    
&$\matrix{dim \Re \cr 5\cr 
dim I \cr 3 \cr }$ && 
$I =\left(\matrix{\alpha&0&0&0\cr 0&\beta&0&0\cr 0&0&
\gamma & 0\cr 0&0&0&\gamma \cr}\right)$
           && 
Perturbation zero &\cr
height5pt&\omit&& \omit&& \omit&\cr
height5pt&\omit&& \omit&& \omit&\cr
\noalign{\hrule }
height2pt&\omit&& \omit&& \omit&\cr
\noalign{\hrule }}}
\section{Acknowledgments}
The author wishes to thank V. Kharchenko for helpful discussion and
CONACYT-M\'exico for partial support under grant No 4336-E.

 \end{document}